\documentclass[proc]{edpsmath}
\usepackage{color}
\usepackage{ulem}
\usepackage{amssymb}

\newtheorem{prop}{Proposition}[section]
\newtheorem{theo}{Theorem}[section]
\newtheorem{lem}{ Lemma}[section]
\newtheorem{corr}{Corollary}[section]

\newtheorem{rem}{Remark}[section]

\def\Var{\mathop{\rm Var}\nolimits}%

\newcommand{\eref}[1]{(\ref{#1})}

\def\1{{{\mbox{${\rm{1\negthinspace\negthinspace I}}$}}}}

\newcommand{\E}{\mathbb{E}}
\newcommand{\PE}{\mathrm{PE}}
\newcommand{\SI}{\mathrm{SI}}
\newcommand{\R}{\mathbb{R}}
\newcommand{\bL}{\mathbb{L}}

\newcommand{\bG}{\mathbb{G}}
\newcommand{\argmin}{\mathrm{argmin}}
\newcommand{\rdg}{\mathrm{rdg}}
\newcommand{\GE}{\mathrm{GE}}
\newcommand{\ER}{\mathrm{ER}}
\newcommand{\cF}{\mathcal{F}}
\newcommand{\T}{\mathrm{GS}}
\newcommand{\cA}{\mathcal{A}}
\newcommand{\cB}{\mathcal{B}}
\newcommand{\cC}{\mathcal{C}}
\newcommand{\cD}{\mathcal{D}}
\newcommand{\cG}{\mathcal{G}}
\newcommand{\cK}{\mathcal{K}}
\newcommand{\cH}{\mathcal{H}}
\newcommand{\cL}{\mathcal{L}}
\newcommand{\cN}{\mathcal{N}}
\newcommand{\cP}{\mathcal{P}}
\newcommand{\cT}{\mathcal{T}}
\newcommand{\cX}{\mathcal{X}}
\newcommand{\gell}{\boldsymbol{\ell}}
\newcommand{\gtheta}{\boldsymbol{\theta}}
\newcommand{\geps}{\boldsymbol{\varepsilon}}
\newcommand{\gbeta}{\boldsymbol{\beta}}
\newcommand{\gX}{\mathbf{X}}
\newcommand{\gx}{\mathbf{x}}
\newcommand{\gy}{\mathbf{y}}
\newcommand{\gY}{\mathbf{Y}}
\newcommand{\gR}{\mathbf{R}}

\newcommand{\PX}{P_{\gX}}
\newcommand{\EX}{E_{\gX}}

\newcommand{\Ee}{E_{\geps}}
\newcommand{\PXe}{P_{\gX, \geps}}
\newcommand{\EXe}{E_{\gX, \geps}}

\begin{document}
\bibliographystyle{plain}
\title{Metamodel construction for Sensitivity Analysis}
\author{Sylvie Huet}
\address{Unité MaIAGE, INRA Jouy-en-Josas, France;
\email{sylvie.huet@inra.fr}}
\author{Marie-Luce Taupin}
\sameaddress{1}
\address{Laboratoire LaMME, UMR CNRS 8071- USC INRA, Universit\'e d'Evry Val d'Essonne, France;\\
\email{marie-luce.taupin@univ-evry.fr}}
\begin{abstract} 
We propose to estimate a metamodel and the sensitivity indices of a
complex model $m$ in the Gaussian regression framework.  Our approach combines methods for sensitivity
analysis of complex models and statistical tools for
sparse non-parametric estimation in multivariate Gaussian regression model. 
It rests on the construction of a
metamodel for aproximating the Hoeffding-Sobol decomposition of $m$. This metamodel belongs to a reproducing
kernel Hilbert space  constructed as a direct sum of Hilbert spaces leading to a functional ANOVA
decomposition. The estimation of the metamodel is carried out  via a penalized
least-squares  minimization allowing  to select the subsets of
variables that contribute to predict the output. It allows to estimate the sensitivity indices of
$m$. We  establish an oracle-type inequality for the risk of the estimator,
describe the procedure for estimating the metamodel and the
sensitivity indices, and assess the performances of the procedure via a
simulation study.
\end{abstract}

\begin{resume} 
Nous consid\'erons l'estimation d'un m\'eta-mod\`ele d'un mod\`ele complexe $m$
\`a partir des observations d'un $n$-\'echantillon dans un mod\`ele de
r\'egression Gaussien. Nous en d\'eduisons une estimation des indices de
sensibilit\'e de $m$. Notre approche combine les m\'ethodes d'analyse de
sensibilit\'e de mod\`eles complexes et les outils statistiques de
l'estimation non-param\'etrique en r\'egression multivari\'ee. Elle repose
sur la construction d'un m\'eta-mod\`ele qui approche la d\'ecomposition de
Hoeffding-Sobol de $m$. Ce m\'eta-mod\`ele appartient \`a un espace de
Hilbert \`a noyau
reproduisant qui est lui-m\^eme la somme directe d'espaces de Hilbert,
permettant ainsi une d\'ecomposition de type ANOVA. On en d\'eduit des
estimateurs des indices de sensibilit\'e de $m$. Nous \'etablissons une
in\'egalit\'e de type oracle pour le risque de l'estimateur, nous
d\'ecrivons la proc\'edure pour estimer le m\'eta-mod\`ele et les indices de
sensibilit\'e, et \'evaluons les performances de notre m\'ethode \`a l'aide
d'une \'etude de simulations.
\end{resume}
\runningtitle{Metamodel construction for Sensitivity Analysis}
\runningauthors{S. Huet and M.L. Taupin}
\maketitle

\section{Introduction}

We consider a Gaussian regression model
\begin{equation}
Y = m(\gX) + \sigma \varepsilon,
\label{RegMod.eq} 
\end{equation}
where $\gX$ is a $d$ random vector with a known distribution $\PX =
P_{1}\times\ldots\times P_{d}$ on $\cX$ a compact subset of $\R^{d}$,
and $\varepsilon$ is independent of $\gX$, and distributed as
a $\cN(0,1)$ variable. The variance $\sigma^{2}$ is unknown and the number of variables $d$ may be large. 
The function $m$  is a complex and unknown function from $\R^{d}$ to
$\R$. It may present strong non-linearities and high order interaction effects
between its coordinates.  On the basis of a $n$-sample $(Y_i, \gX_i), i=1,
\ldots, n$, we aim to  construct metamodels
and perform sensitivity analysis in order to determine
the infuence of each variable or group of variables on the output. 
 
Our approach combines methods for sensitivity
analysis of a complex model and statistical tools for
sparse non-parametric estimation in multivariate Gaussian regression model. 
It rests on the construction of a
metamodel for aproximating the Hoeffding-Sobol decomposition of the
function $m$. This metamodel belongs to a reproducing
kernel Hilbert space  constructed as a direct sum of Hilbert spaces leading to a functional ANOVA
decomposition involving variables and interactions between them. The estimation of the metamodel is carried out  via a penalized
least-squares  minimization allowing  to select the subsets of
variables $\gX$ that contribute to predict the output $Y$. Finally, the
estimated metamodel allows to estimate the sensitivity indices of $m$.

A lot of work has been done around meta-modelling and sensitivity
indices estimation.

For a complete account on global sensitivity analysis, see for example the book
by Saltelli et al.~\cite{Saltellietal}. Let us briefly present the context of usual global sensitivity analysis. Suppose that we are able
to calculate the ouputs $z$ of a  model $m$ for $n$ realizations of
the input vector  $\gX$, such that $z_i = m(\gX_i)$ for $i=1, \ldots,
n$. Starting from the values  $(z_i, \gX_i),
i=1, \ldots n$,  the objectives of meta-modelling and global
sensitivity analysis are to approximate the function $m$ by what is
called a metamode or are to quantify the influence of 
 some subsets of
the variables $\gX$ on the output $z$. This metamodel  helps to understand the behavior of the
model, or allows to speed up future calculation using it in place
of the original model $m$.  

In particular when
the inputs variables $\gX$ are  independent, if $m$
is  square integrable, one may
consider the classical Hoeffding-Sobol decomposition~\cite{Sobol1993,VdVaart98} that leads to write
$m$ according to  its ANOVA functional expansion:
\begin{equation}
m(\gx) = m_{0} + \sum_{v \in \cP} m_{v} (\gx_{v}),
\label{hoeffding.eq}
\end{equation}
where  $\cP$ denotes the set of parts of
$\left\{1, \ldots, d\right\}$ with dimension 1 to $d$ and where for
all $\gx \in \R^{d}$, $\gx_{v}$ denotes the vector with components
$x_{j}$ for $j \in v$.  The functions $m_v$ are centered and orthogonal
in $\bL^{2}(\PX)$ leading to the following decomposition of the
variance of $m$: $\Var\left(m(\gx)\right) =
\sum_{v}\Var\left(m_{v}(\gx_{v})\right)$. The Sobol 
sensitivity indices,  introduced by Sobol~\cite{Sobol1993}, are
defined for any group of variables $\gx_v$, $v\in 
\cP$ by 
\begin{equation*}
 S_{v} = \frac{\Var\left(m_{v}(\gx_{v})\right) }{\Var\left(m(\gx)\right)}.
\end{equation*}
They quantify
the contribution of a subset of variables $\gx$ to the output
$m(\gx)$.
Several approaches are available for estimating these
sensitivity indices, see for example Iooss and
Lema\^itre~\cite{ioossLemaitre15} for a recent review. Among all of
them, let us consider  the one based on metamodel construction
that allows to directly obtain the sensitivity indices. 
Generally one consider metamodels that correspond to an ANOVA-type
decomposition, and that are candidate to approximate the Hoeffding
decomposition of $m$. 
The ANOVA-type decomposition leads to
consider functions defined as follows: 
\begin{equation}
f : \cX \rightarrow \R, f( \gx) = f_{0} + \sum_{v \in \cP}
f_{v}(\gx_{v}),
\; \E_{\PX} f_{v}(\gx_{v}) = \E_{\PX} f_{v}(\gx_{v} )f_{v'}(\gx_{v'} )=0
\; \forall v, v' \in \cP
\label{anova.eq}
\end{equation}
for functions $f_{v}$ that are chosen to belong to some
functionnal spaces. The  polynomial
Chaos construction, see for example Ghanem and
Spanos~\cite{GhanemSpanos}, Soize and Ghanem~\cite{SoizeGhanem2004}, can be used to approximate the Hoeffding
decomposition of $m$. 
This approach  was considered by
Blatman and Sudret~\cite{BlatmanSudret} who propose a method
for truncating (such that to keep
polynomials with degree less 
than some integer)  the polynomial Chaos expansion and then an algorithm
based on least angle regression for selecting the terms in the
expansion. For the same purpose, Gu and Wu~\cite{GuWu} propose an
algorithm based on the hierarchy principle (lower
order effects are more likely to be important than higher order
effects) and on the  heredity principle (interaction can be active only if one or all of
its parent effects are also active). This approach joins the one
proposed by Bach~\cite{bach09} for variable selection based on
hierarchical kernel learning.

Inspired by Touzani \cite{Touzani}, Durrande et al.~\cite{Durrande2013} propose to approximate
 $m$ by functions belonging to a reproducing kernel Hilbert space
 (RKHS). The RKHS is constructed as a direct sum of Hilbert spaces
 leading to a functional ANOVA decomposition (see Equation~\eref{anova.eq}), such that the
 projection of $m$   onto the RKHS is an approximation of the
 Hoeffding decomposition of $m$. 

Following Lin and Zhang~\cite{LinZhang06}, Touzani and Busby~\cite{Touzani13} propose an algorithm to calculate the
 penalized least-square estimator of $m$ on the RKHS space, where the
 least-square criteria is penalized by the sum 
 of the norms  of $m$ on each Hilbert subspace. This group-lasso type procedure
 allows both to select and calculate the non-zero terms in the functional ANOVA
 decomposition. 

\medskip
\noindent
Our objective is to propose an estimator of a metamodel which will
approximate the Hoeffding decomposition of $m$ considering a  Gaussian
regression model defined at Equation~\eref{RegMod.eq} and to deduce
from this estimated metamodel, estimators for the sensitivity
indices of $m$. Contrary to the usual setting of sensitivity analysis where $m(\gX_i)$ is available, only the observations $Y$ are
available, which leads us to consider  the nonparametric multivariate regression setting.

Let us briefly describe  the methods related to this regression setting and review their
theoretical properties, starting with papers assuming an  univariate additive decomposition for the  function  $m$ in the context of high-dimensional
sparse additive models. Precisely,  denote by  $\cF^{\textrm{1-add}}$, the set of functions $f$ defined on
$\cX$ such that 
$f(\gx) =f_{0} + \sum_{a=1}^{d} f_{a}(x_{a})$ 
where $f_{0}$ is a constant, and where for $a=1, \ldots, d$, the  functions $f_{a}$ are centered
and square-integrable with respect to $P_{a}$. For
each function $f$, the set $S_{f}$ of  
indices $a \in \left\{1, \ldots d\right\}$ such that $f_{a}$ is
not identically zero is called the actice set of $f$.

Ravikumar et al.~\cite{Ravikumar09} 
propose a group-lasso procedure where each function $f_a$ is
approximated by its truncated decomposition on a basis of
functions. They provide an algorithm and, assuming that the function
$m$ belongs to the set $\cF^{\textrm{1-add}}$ and that $S_{m}$  is sparse, prove
the consistency  of the 
active set and of the risk  of the estimator of $m$.

Meier et al.~\cite{meier2009} propose to combine a sparsity
penalty (group-lasso) and a smoothness penalty (ridge) for estimating
$m(\gx)$. Considering the fixed design
framework, they 
establish some oracle properties of the empirical risk for estimating
the projection of $m$ onto the set of univariate additive functions
$\cF^{\textrm{1-add}}$. 

Raskutti et al.~\cite{RWY} consider the case where each univariate
function $f_{a}$ belongs to a RKHS and as Meier et al. combine a sparsity and a smoothness
penalty. 
Assuming that the $d$ variables $\gX$ are independent, they
derive upper bounds for the integrated and the empirical risks, as
well as a lower bound for the integrated risk over spaces of sparse
additive models whose each component is bounded with respect
to the RKHS norm. 

Additive sparse modelling is  too restrictive
in practical settings because it  does not take into account
interactions between variables that may affect
the relationship between $Y$ and $\gX$. The generalization of additive
smoothing splines to interaction smoothing splines leading to an
ANOVA-type decomposition (see Equation~\eref{anova.eq})  was proposed by several
authors (see for example Wahba~\cite{Wahba90},
Friedman~\cite{Friedman91}, Wahba et
al.~\cite{Wahbaetal95}).

To control smoothness and to enforce
sparsity in the ANOVA-type decomposition, Gunn and
Kandola~\cite{GunnKandola} propose to consider the ANOVA
decomposition as a weighted linear sum of kernels and to use a lasso
penalty on the weights to select the terms in the decomposition as
well as a ridge penalty to ensure smoothness of the kernel expansion.
The COSSO proposed by Lin and Zhang~\cite{LinZhang06}
is based on  smoothness penalty defined as the sum of the RKHS-norms of the
functions $f_{v}$. The authors study existence and rate of convergence
of the estimator. In a more general framework, where the function $m$
 is written as a
linear span of a large number of kernels, Koltchinskii and
Yuan~\cite{koltchinskiiYuan08} established oracle inequalities on
the excess risk assuming that the function $m$ has a sparse
representation (the set of $v \in \cP$ such that $f^{*}_{v}$ is non
zero is sparse). The  authors generalized their results
to a penalty function that combines sparsity and smoothness~\cite{koltchinskii2010}, as
proposed by Meier et al.~\cite{meier2009} and Raskutti et
al.~\cite{RWY}. 

Recently Kandasamy and Wu~\cite{KandasamyYu2016} proposed an
estimator called SALSA, based on a ridge penalty, where the
ANOVA-type decomposition is restricted to set $v \in \cP$ such
that $|v| \leq D_{\max}$. The authors propose to choose $D_{\max}$ via
a cross-validation procedure.

\subsection*{Our contributions}

Using the
functionnal ANOVA-type decomposition as proposed by Durrande et
al.~\cite{Durrande2013}, we propose an estimator of a metamodel which 
approximates the Hoeffding decomposition of $m$. 
Following 
the most recent works  in the framework of nonparametric estimation of
sparse additive models, we propose a penalized least-square estimator
where the penalty function enforces both the sparsity and the
smoothness of the terms in
the decomposition. We show that our estimator satisfies an oracle
inequality with respect to the empirical and integrated risks. 

Our procedure allows both to select and estimate  the terms in the  ANOVA decomposition, and therefore,  to select the sensitivity indices that are non-zero and
estimate them. In particular it makes possible to estimate  Sobol indices of high order, a point  known to be difficult in practice.

Finally, using convex optimization tools,  we develop an algorithm (available on request) in $R$ \cite{R},   for calculating the estimator. A simulation study 
shows the good performances of our estimator in practice.
 
The paper is organised as follows: The RKHS construction based on ANOVA
kernels and the procedure for estimating a
metamodel are presented in Section~\ref{procedure.st}. The estimators
of the Sobol indices are given in Section~\ref{sensA.st}.  The theoretical 
properties of the metamodel estimator are stated in 
Theorem~\ref{oracle} and Corollaries~\ref{oracle2} and~\ref{oracle3} whose
proofs are postponed in Sections~\ref{sketch.st} and~\ref{proofs.st}.  Section~\ref{CalcEstim.st} is devoted to the
calculation of the estimator and  Section~\ref{simul.st} to the
simulation study.


\section{Meta-modelling}
We start from the Hoeffding decomposition (see Sobol~\cite{Sobol2001} and
Van der Vaart~\cite{VdVaart98}, p. 157) of the 
function $m$ that consists in writting $m$ as in Equation~\eref{hoeffding.eq}
\begin{equation*}
m(\gx) = m_{0} + \sum_{v \in \cP} m_{v} (\gx_{v}),
\end{equation*}
where  $\cP$ denotes the set of parts of
$\left\{1, \ldots, d\right\}$ with dimension 1 to $d$ and where for
all $\gx \in \R^{d}$, $\gx_{v}$ denotes the vector with components
$x_{j}$ for $j \in v$.  
For all $v, v'$ in $\cP$, 
\begin{equation*}
 \EX\left(m_{v}(\gX_{v})\right) =
\EX\left(m_{v}(\gX_{v}) m_{v'}(\gX_{v'})\right) = 0.
\end{equation*}

\label{procedure.st}
We propose to consider a functionnal space based on the tensorial
product of Reproducing Kernel Hilbert spaces (RKHS), and to
approximate the unknown function $m$ by its projection denoted $f^*$
on such such RKHS space. 
One of the key point is to construct the space $\cH$ such that the
terms of the decomposition of a function $f$ in $\cH$ correspond to
its Hoeffding-Sobol decomposition.

\subsection{RKHS construction}
\label{construction.st}
Let us describe the construction of spaces $\cH$, based on ANOVA kernels, construction which was
given by Durrande et 
al.~\cite{Durrande2013}.

Let $\cX = \cX_{1} \times \ldots \times \cX_{d}$ be a compact subset
of $\R^{d}$.
For each coordinate $a\in \{1,\cdots,d\}$, we choose a RKHS  $\cH_{a}$ and
its associated kernel $k_{a}$ defined on the set $\cX_{a} \subset \R$ such
that
the two following properties are satisfied
\begin{enumerate}
\item $k_{a} : \cX_{a} \times \cX_{a} \rightarrow \R$ is $P_{a} \times
  P_{a}$ measurable
\item $\E_{P_{a}}\sqrt{k_{a}(X_{a}, X_{a})} < \infty$
\end{enumerate}
The RKHS $\cH_{a}$ may be decomposed as $\cH_{a} = \cH_{0 a}
\stackrel{\perp}{\oplus} \cH_{1 a}$, where 
\begin{eqnarray*}
 \cH_{0 a} & = & \left\{ f_{a} \in \cH_{a}, \E_{P_{a}}(f_{a}(X_{a})) =
   0\right\}\\
\cH_{1 a} & = & \left\{  f_{a} \in \cH_{a}, f_{a}(X_{a}) =C \right\},
\end{eqnarray*}
the kernel $k_{0a}$ associated to the RKHS $\cH_{0 a}$ being defined
as follows (see Berlinet et Thomas-Agnan~\cite{berlinet2003reproducing}):
\begin{equation*}
k_{0a} (x_{a},x'_{a}) = k_{a}(x_{a},x'_{a}) - 
\frac{\E_{U \sim P_{a}}(k_{a}(x_{a},U))\E_{U \sim P_{a}}(k_{a}(x'_{a},U))}
{\E_{(U,V)\sim P_{a}\times P_{a}}k_{a}(U,V)}.
\end{equation*}
The  ANOVA kernel is defined as 
\begin{equation*}
 k(\gx, \gx') = \prod_{a=1}^{d} 
\left(1+k_{0a}(\gx_{a}, \gx'_{a})\right) = 
1 + \sum_{v \in \cP} k_{v}(\gx_{v}, \gx'_{v}), \mbox{ where }
k_{v}(\gx_{v}, \gx'_{v}) = \prod_{a \in v} k_{0a} (x_{a},x'_{a}),
\end{equation*}
and the corresponding RKHS 
\begin{equation*}
 \cH = \otimes_{a=1}^{d}\left( 1 \stackrel{\perp}{\oplus}
   \cH_{0a}\right) = 1 + \sum_{v \in \cP} \cH_{v},
\end{equation*}
where the RKHS $\cH_{v}$ is associated with kernel $k_{v}$.
According to this construction, any function $f \in \cH$ satisfies
\begin{equation*}
f(\gx) = \langle f, k( \gx, \cdot) \rangle_{\cH} = 
f_{0} + \sum_{v \in \cP} f_{v}(\gx), 
\end{equation*}
where $f_{v}(\gx) = \langle f, k_{v}( \gx, \cdot) \rangle_{\cH}$
depends on $\gx_{v}$ only. For all $v \in \cP$, $f_{v}(\gx_{v})$ is 
centered and for all $v' \neq v$, $f_{v}(\gx_{v})$ and
$f_{v'}(\gx_{v'})$ are uncorrelated. We get thus the Hoeffding
decomposition of $f$.

\subsection{Approximating the Hoeffding decomposition of $m$}
\label{approxHoeffding.st}
  
Let
$f^{*}=f_{0}^{*} + \sum_{v\in \cP} f_{v}^{*}$ which minimizes 
\begin{equation*}
 \|m-f\|^{2}_{\bL^{2}(\PX)} = \EX\left(m(\gX)-f(\gX)\right)^{2}
\end{equation*}
over functions $f \in \cH$. This $f^*$ can be viewed as an approximation of $m$
and more specifically his Hoeffding decomposition is an approximation of
the Hoeffding decomposition of
$m$. Therefore if the Hoeffding decomposition of $m$ is written as in
Equation~\eref{hoeffding.eq}, each function $f^{*}_{v}$ approximates the function $m_{v}$. 

The idea is propose an estimator of $f^*$ as estimator of $m$.

\subsection{Selection step}
\label{select.st}

Since $\cP$ is the set of parts of
$\left\{1, \ldots, d\right\}$, the number of functions$f^{*}_{v}$ is related to the cardinality of $\mathcal{P}=2^d-1$ that may be huge.  Our construction is thus associated to a selection strategy.

The selection of  $f^{*}_{v}$ in $f^{*}$ is based on a {\it ridge-group-sparse} type
procedure which minimizes the  penalized least-squares criteria over
functions $f \in \cH$. The least-squares criteria is penalized in order to
both select few terms in the additive decomposition of $f$ over sets $v \in
\cP$, and to favour smoothness of the estimated $f_{v}$. The ridge
regularization is  ensured by controling the norm of $f_{v}$ in the
Hilbert space $\cH_{v}$ for all $v$, and the group-sparse
regularization is strengthened by controling the empirical norm of
$f_{v}$, defined as 
\begin{equation*}
 \|f\|_{n} = \sqrt{\frac{1}{n} \sum_{i=1}^{n} f^{2}_{v}(\gX_{v, i})}.
\end{equation*}
For any $f \in \cH$ such that $f=f_{0} + \sum_{v\in \cP } f_{v}$, and for some tuning parameters $(\mu_{v}, \gamma_{v}, v \in \cP)$, let 
$\cL(f)$ be defined as 
\begin{equation}
\cL(f) = \frac{1}{n} \sum_{i=1}^{n} 
\left(Y_{i} - f_{0} - \sum_{v \in \cP} f_{v}(\gX_{v, i})\right)^{2}
+ \sum_{v \in \cP} \mu_{v} \|f_{v}\|_{\cH_{v}} + 
\sum_{v \in \cP} \gamma_{v}\|f_{v}\|_n.
\label{PenCrit.eq}
\end{equation}
Let us define the set of functions $\cF$
\begin{equation}
 \cF = \left\{ f \mbox{ such that } f=f_{0} + \sum_{v\in \cP} f_{v},
   f_v \in \cH_{v}, \|f_v\|_{\cH_v} \leq 1 \right\}.
\label{calF.eq}
\end{equation}
Then
the estimator $\widehat{f}$ is defined as
\begin{equation}
 \label{fchap}
\widehat{f} = \argmin \left\{ \cL(f), f \in \cF \right\}.
\end{equation}

\begin{rem}
The construction of the RKHS spaces  described
above, allows to consider functionnal 
spaces that suit well to the smoothness of the function $m$,
irrespectively of the distribution of $\gX$. Indeed, the kernels
$k_{0,a}$ depend on the distribution of $\gX$ only for calculating the projection
onto the space of constant functions. 
In comparison, the decomposition based on
the truncated polynomial Chaos expansion, used for
sensitivity analysis (see Blatman and Sudret~\cite{BlatmanSudret}), is
based on the distribution 
of $\gX$, and only the choice of the truncation handles the smoothness
of the approximation. 
\end{rem}

\section{Sensitivity analysis }
\label{sensA.st}

\subsection{Sobol indices}
\label{sobol.st}


Let us go back to the  Hoeffding decomposition Equation~\eref{hoeffding.eq}.
The orthogonality between two terms in this decomposition
leads to the additive decomposition of the
variance of $m(\gx)$:
\begin{equation*}
\Var\left(m(\gx)\right) = \sum_{v \in \cP} \Var\left(m_{v}(\gx_{v})\right).
\end{equation*}
Each of these variance terms are related to Sobol indices~\cite{Sobol1993}. For
example, the Sobol indice linked with the interaction between
variables $\gx_{v}$ is defined as 
\begin{equation*}
 S_{v} = \frac{\Var\left(m_{v}(\gx_{v})\right) }{\Var\left(m(\gx)\right)},
\end{equation*}
or the global Sobol indices for the variable $x_{a}$, $a\in \{1,\cdots,d\}$, is
\begin{equation*}
G_{a}  = \frac{\sum_{v\supseteq \left\{a\right\}}\Var\left(m_{v}(\gx_{v})\right) }{\Var\left(m(\gx)\right)}.
\end{equation*}
Those Sobol indices and global Sobol indices quantify
the contribution of a subset of variables $\gx$ to the output
$m(\gx)$.
As it is said in the introduction direct estimation of these Sobol indices may require lot of calculations. We consider here methods based on metamodels to directly obtain Sensitivity indices.

\subsection{Estimation of Sobol indices }
\label{estsobol.st}

Thanks to the orthogonality property of functions in $\cH$, the variance of $m(\gx)$  will be estimated by 
\begin{equation}
 \widehat{\Var}\left(m(\gx)\right) = \sum_{v \in \cP}
 \widehat{\Var}\left(m_{v}(\gx_{v})\right), \mbox{ where } 
\widehat{\Var}\left(m_{v}(\gx_{v})\right) =
\E_{\gX}\left(\widehat{f}^{2}_{v}(\gX_{v})\right) = \|
\widehat{f}_{v}\|^{2}_{\bL^{2}(\PX)}.
\label{VarChap.eq}
\end{equation}

In practice, in order to avoid calculating the variance of
$\widehat{f}_{v}(\gX_{v})$, one may use an  estimator based on the empirical variances of functions
$\widehat{f}_{v}$. Precisely, if $\widehat{f}_{v,\cdot}$ is the mean of  the
 $\widehat{f}_{v}(\gX_{v,i})$, for $i=1, \ldots, n$, then 
\begin{equation}
\widehat{\Var}^{{\rm emp}}\left(m_{v}(\gx_{v})\right) =  \frac{1}{n-1}
\sum_{i=1}^{n} \left( 
\widehat{f}_{v}(\gX_{v,i}) - \widehat{f}_{v, \cdot}\right)^{2}
\label{Vemp.eq}.
\end{equation}

One of the main contribution of this approach is to allow the
estimation of Sobol indices of any order, whereas classical methods
only deal with small order, generally less than  two.

\section{\label{ThRes.st}Theoretical result: oracle inequality for metamodel}

In this section we state 
 the oracle inequality  for the estimated metamodel $\widehat{f}$ which approximates the Hoeffding decomposition of the
  unknown function $m$.  
  
  \subsection{Notations and Assumptions}
  \label{notation.st}
  
For a function $f \in \cH$, $f=f_{0} + \sum_{v\in \cP} f_{v}$, we
denote by $S_{f}$ its support and $|S_{f}|$ its cardinality. More precisely
\begin{equation}
S_{f}  = \left\{ v \in \cP, f_v \neq 0 \right\}.
\label{sparsity}
\end{equation}

We consider RKHS spaces $\cH_{v}, v \in \cP$ satisfying the following assumptions:
\begin{itemize}
\item for all $f_v \in \cH_v$, $\EX f^{2}_{v}(\gX) < \infty$, 
\item for all $f_v \in \cH_v$, $f_{v'} \in \cH_{v'}$, $\EX f_{v}(\gX)
  = 0$ and  $\EX f_{v}(\gX) f_{v'}(\gX) = 0$,
\item there exists $R' >0$ such that 
\begin{equation}
 \forall f_v \in \cH_v \;\;
\|  f_v\|_{\infty} = \sup\left\{ |f_v(\gX)|, \gX \in \cX \right\} \leq
  R'. \label{supfv.eq}
\end{equation}
\end{itemize}

For each $v \in \cP$, let  $\omega_{v,k}$, for $ k\geq 1$ be the eigenvalues of the operator
associated to the self reproducing kernel $k_v$, arranged in the
decreasing order. Let us define the function $Q_{n,v}(t)$, for positive
$t$, as follows:
\begin{equation}
\label{Qn}
Q_{n,v}(t)=\sqrt{\frac{5}{n} \sum_{k\geq 1} \min(t^2,  \omega_{v,k}  )},
\end{equation}
and for some $\Delta > 0$  let $\nu_{n,v}$ 
be defined as follows
\begin{eqnarray}
\label{nu}
\nu_{n, v} = \inf\left\{ t \mbox{ such that } Q_{n,v}(t)\leq  \Delta t^{2}\right\}.
\end{eqnarray}
For each $v\in \mathcal{P}$, $\nu_{n,v}$ refers to the so-called critical univariate rate, the minimax-optimal rate for $\mathbb{L}^2(P_\gX)$-estimation of a single univariate function in the hilbert space $\mathcal{H}_v$ (e.g. Mendelson~\cite{Mendelson2002}).

Our choices of regularization parameters and rates are specified in terms of the quantities:
\begin{eqnarray}
 \lambda_{n,v} &= &\max\left\{
  \nu_{n,v}, \sqrt{d/n} \right\} \label{lambda.eq}.
\end{eqnarray}

\begin{theo}
\label{oracle}
Let us consider the regression model defined at Equation~\eref{RegMod.eq}.
Let  $(Y_{i}, \gX_{i})$, $i=1, \ldots, n$ be  a $n$-sample  with
the same law as $(Y, \gX)$. Let $\widehat{f}$ be defined
by~\eref{fchap}. 

If there exist
constants $C_{l}, l=1, 2, 3$, $C_{1} \geq 1$, and $0 < \eta <1$ such that
 the following conditions are satisfied:
\begin{equation}
\label{condmu}
\mbox{for all } v \in \cP, 
\lambda_{n,v} \leq \min\left\{ 
\frac{\gamma_{v}}{C_{1}}, \sqrt{\frac{\mu_{v}}{C_{1}}}\right\}, \;
n\lambda_{n,v}^{2} \geq -C_{2} \log\lambda_{n,v}, \;
\end{equation}
and 
\begin{equation}
\mbox{ for all } f \in \cF,  \; \max\left( \sum_{v\in
    S_f} \gamma_v^2,\sum_{v\in S_f}\mu_v   \right) \leq 1,
\label{hypdf.eq}
\end{equation}
then, with probability greater than $1-\eta$,
\begin{equation*}
 \| m - \widehat{f}\|^{2}_{n} \leq
C_{3} \inf_{f \in \cF}\left\{  \| m - f\|^{2}_{n}
+ \sigma^{2} \sum_{v \in S_{f}} \left(\mu_{v} + \gamma_{v}^{2}\right)
 \right\}.
\end{equation*}

\end{theo}

\medskip
\noindent
The result can be easily generalized to the minimization of $\cL(f)$
over the space $\cG$ defined as follows:
For some positive constants $r_{v}, v \in \cP$, 
\begin{equation}
 \cF_r = \left\{ f \mbox{ such that } f=f_{0} + \sum_{v\in \cP} f_{v},
   f_v \in \cH_{v}, \|f_v\|_{\cH_v} \leq  r_{v} \right\}.
\label{Fr.eq}
\end{equation}
Indeed, we just have to consider the RKHS $\cH^{'}_{v}$ associated
with the kernel $k^{'}_{v}(\gx_{v},\gx^{'}_{v} ) = r_{v}^{2}
k_{v}(\gx_{v},\gx^{'}_{v} )$ in place of the RKHS $\cH_{v}$. Then minimizing 
\begin{equation*}
 \frac{1}{n} \sum_{i} \left(
Y_{i} - f_{0} - \sum_{v} f_{v}(X_{v,i})\right)^{2}
+ \sum_{v} \mu_{v}  r_{v} \|f_{v}\|_{\cH^{'}_{v}} 
+ \sum_{v} \gamma_{v} \|f_{v}\|_{n}
\end{equation*}
over the set
\begin{equation*}
\cF^{'} = \left\{ f \mbox{ such that } f=f_{0} + \sum_{v} f_{v}, f_{v} \in \cH^{'}_{v}, \|f_{v}\|_{\cH^{'}_{v}}
   \leq 1
\right\}.
\end{equation*}
is equivalent to minimizing $\cL(f)$ over the set $\cF_r$.

\medskip
\noindent
Let us now comment on the theorem.
\begin{itemize}
\item The term $\| m -f\|_n^2$ refers to the usual bias term quantifying the approximation properties of the Hilbert space  $\mathcal{H}$ as a distance between
the true $m$ and $f$, its approximation into $\mathcal{H}$.

\item Koltchinskii et Yuan~\cite{koltchinskii2010} considered what is
  called the multiple kernel learning problem, where the functions in
  $\cH$ have an additive representation over kernel spaces. They do
  not assume that the variable $\gX$ are independent, nor that the
  kernel spaces satisfy an orthogonality condition. In return, they
  assume that some decomposability type properties are satisfied, and
  they introduce some characteristics related to the degree of
  ``dependence'' of the kernel spaces. 

\item In the particular case when the decomposition is limited to the
  main effects of the variables, then the problem comes back to the
  classical nonparametric additive model. The theoretical properties
  of the estimator based on a {\it ridge-group-sparse} type procedure
  have already  been established  (see for example Meier et
  al.~\cite{meier2009}, and 
    Raskutti et al.~\cite{RWY}).

\item Weights in the penalty terms may be of interest for
  applications. The theoretical result highlights that the tuning parameters $(\mu_{v},
  \gamma_{v})$ should depend on
  the decreasing of the eigenvalues of the kernel defining
  $\cH_{v}$. Besides, we may be interested by introducing weigths that
  favor small order interaction terms. 
\item Because we aim to approximate the Hoeffding
  decomposition of $m$, we need to have 
  orthogonality between the spaces $\cH_v, v \in \cP$. This condition,
  required by our objective, is also  a
  key point in the proof of the Theorem, when the problem is to
  compare the euclidean norm of functions in $\cH$ with the
  norm in $\bL^{2}(\PX)$. At this step we need to assume that
  Assumption~\eref{hypdf.eq} holds to conclude.

\item 
We do not
  assume that the functions in $\cF$ are uniformly globally bounded,
  that is that the $\sup
  \{ |f(\gx)|,\gx \in \cX \}$ is bounded by a constant that does not
  depend on $f$. 
  Instead we assume that each function within the unit ball of the Hilbert space $\mathcal{H}_v$ is uniformly bounded by a constant multiple of its Hilbert norm. 
 In fact, the functions $f$ in the
  space $\cF$, written as $f_{0} + \sum_{v} f_{v}$ satisfy that for
  each group $v$, $f_{v}$ is uniformly bounded. This assumption  is easily satisfied as soon as the kernel $k_v$ is bounded
  on the compact set $\cX$.  Indeed, $\| f_v\|_{\infty} \leq
  \sup_{\gX \in \cX} \sqrt{k_{v}(\gX_v, \gX_v)} \|f_v\|_{\cH_{v}}$. We refer to~\cite{RWY} for
  a discussion on that subject and a comparison with the work of
  Koltchinskii et Yuan~\cite{koltchinskii2010}.

\item Assuming that $n\lambda_{n,v}^{2} \geq -C_{2} \log\lambda_{n,v}$
  allows to control the probability of the union of $|\cP|$
  events. This is a mild condition, satisfied for $\lambda_{n,v} =
  K_{n}/\sqrt{n}$ for  $K_{n}$ of the order $\sqrt{\log n}$.
\end{itemize}

\medskip
\noindent
The following corollary  gives an upper bound of the risk with respect
to the 
$\bL^{2}(\PX)$ norm. It is  mainly a consequence of
Theorem~\ref{oracle} (its proof is given Section~\ref{PrCorr2.1.st} page~\pageref{PrCorr2.1.st}).

\begin{corr}
\label{oracle2}
Under the assumptions of Theorem~\ref{oracle}, we have that, with
probability greater than $1-\eta$, for some constant $C_{4}$,
\begin{equation*}
 \| m - \widehat{f}\|_{\bL^{2}(\PX)}^{2} \leq
C_{4} \inf_{f \in \cF}\left\{  \| m - f\|^{2}_{n} +  \| m -  f\|_{\bL^{2}(\PX)}^{2}
+ \sigma^{2} \sum_{v \in S_{f}} \left(\mu_{v} + \gamma_{v}^{2}\right)\right\}.
\end{equation*}
\end{corr}

From this corollary we can compare the
$\widehat{\Var}\left(m_{v}(\gx_v)\right)$ (see Equation~\eref{VarChap.eq}) with the
variance of $m_{v}(\gx_v)$. Thanks to the following inequality
\begin{equation*}
\left|   \|\widehat{f}_{v}\|_{\bL^{2}(\PX)} - \|m_{v}\|_{\bL^{2}(\PX)}
\right|
\leq  \|\widehat{f}_{v}-m_{v}\|_{\bL^{2}(\PX)} \leq \|\widehat{f}-m\|_{\bL^{2}(\PX)},
\end{equation*}
and to Corollary~\ref{oracle2}, we get the following result
\begin{eqnarray*}
\mbox{if } m_v \equiv 0 & \mbox{then} &
 \widehat{\Var}\left(m_v(\gx_v)\right) \leq
 \|\widehat{f}-m\|^{2}_{\bL^{2}(\PX)} \\
\mbox{if } \|m_v(x_v)\|_{\bL^{2}(\PX)} \geq c >0& \mbox{then} &
\left|\frac{\widehat{\Var}\left(m_v(\gx_v)\right)
}{\Var\left(m_v(\gx_v)\right)} -1\right| \leq \|\widehat{f}-m\|^{2}_{\bL^{2}(\PX)}/c.
\end{eqnarray*}

\subsection{Rate of convergence}
\label{RateCvge.st}

\begin{corr}
\label{oracle3}
Under the same condition as Theorem~\ref{oracle}, if $\gamma_{v} = c
\lambda_{n,v}$ and $\mu_{v} = c \lambda_{n,v}^{2}$ for $c \geq C1$,
then 
\begin{equation*}
 \| m - \widehat{f}\|^{2}_{n} \leq
C_{3} \inf_{f \in \cF}\left\{  \| m - f\|^{2}_{n}
+ \left(\sum_{v \in S_{f}} \nu^{2}_{n,v} + \frac{d |S_{f}|}{n}\right)
\sigma^{2} \right\}.
\end{equation*}
\end{corr}

The result is non asymptotic in the sense that it is shown for any $(n, d)$. Nevertheless, the upper
  bound is relevant when the infimum is reached for functions $f$
  whose decomposition in $\cH$ is sparse, and when $d$ is small face
  to $n$. In fact, the coefficient $d$ occuring in the rate $d
  |S_{f}|/n$ comes from the logarithm of the cardinality of
  $\cP$ equal to $\log(2^{d}-1)$. When $d$ is large, it may be
  judicious to limit the 
  decomposition of functions in $\cH$, to interactions of limited
  order, so that the number 
  of terms in the decomposition stays of the order $\log(d)$.

Let us  discuss the rate of convergence given by
$\sum_{v \in S_{f}} \nu^{2}_{n, v}$. 
For the sake of simplicity let us consider the case where the
variables $X_{1}, \ldots, X_{d}$ have the same distribution $P_{1}$
on $\cX_{1} \subset\R$, and
where the 
unidimensionnal kernels $k_{0 a}$ are all identical, such that
$k_{v}(\gx_{v}, \gx'_{v}) = \prod_{a \in v}
k_{0}(x_{a},x'_{a})$. The kernel $k_{0}$ admits an eigen expansion
given   by 
\begin{equation*}
 k_{0}(x,x') = \sum_{\ell \geq 1} \omega_{0, \ell} \zeta_{\ell}(x)\zeta_{\ell}(x') 
\end{equation*}
where the eigenvalues $\omega_{0, \ell}$ are non negative and ranged
in the decreasing order, and where the $\zeta_{\ell}$ are the
associated eigen functions, orthonormal with respect to $\bL^{2}(P_1)$.
Therefore the kernel $k_{v}$  admits the following expansion
\begin{equation*}
 k_{v}(\gx_v,\gx^{'}_{v}) = \sum_{\gell=(\ell_1 \ldots \ell_{|v|})} 
\underbrace{\prod_{a=1}^{|v|}\omega_{0, \ell_{a}}}_{\omega_{v, \gell}}
\underbrace{\prod_{a=1}^{|v|}\zeta_{\ell_{a}}(x_{a})}_{\zeta_{v, \gell}(\gx_{v})}
\underbrace{\prod_{a=1}^{|v|}\zeta_{\ell_{a}}(x^{'}_{a})}_{\zeta_{v, \gell}(\gx^{'}_{v})}.
\end{equation*}
Consider now the case where the eigenvalues $\omega_{0, \ell}$ are
decreasing at a rate $\ell^{-2\alpha}$ for some $\alpha > 1/2$. 
It can be shown, see Section~\ref{ratecvg.proof}, that the rate
$\nu_{n, v}$ defined at Equation~\eref{nu} is bounded above by a term of
order $n^{-\alpha/(2\alpha + 1)} (\log n)^{\gamma}$, where $\gamma \geq
(|v|-1)\alpha/(2\alpha-1)$. Note that in this particular case, the
rate of convergence depends on $|v|$ through the logarithmic term, and
that up to this logarihmic term the rate of convergence has the same
order than the usual nonparametric  rate for unidimensionnal functions.
 It follows that the RKHS space
$\cH$ should be chosen such that the unknown function $m$ is well
approximated by sparse functions in $\cH$ with low order of
interactions.

\section{Calculation of the estimator}
\label{CalcEstim.st}
The functional minimization problem described at
Equation~\eref{fchap} is equivalent to a parametric minimization
problem. Indeed, we know that if $\cH$ is a RKHS associated with a
kernel $k : \cX \times \cX \rightarrow \R$, then for all $(\gx_{1}, \ldots,  \gx_{n}) \in
\cX^{n}$, and for all $(\alpha_{1}, \ldots, \alpha_{n}) \in \R^{n}$,
the function $f(\cdot) = \sum_{i=1}^{n} \alpha_{i} k(\gx_{i},\cdot)$ is in
$\cH$ and $\|f\|^{2}_{\cH} = \sum_{i,i'=1}^{n} \alpha_{i}\alpha_{i'}
k(\gx_{i},\gx'_{i'})$. In particular, it can be shown that the solution to our
minimization problem is written as $f=f_{0}+\sum_{v\in \cP} f_{v}$
where, according to the representer Theorem (see Kimeld{\?o}rf and Whahba \cite{kimeldorf1970}), $f_{v}(\cdot) = \sum_{i=1}^{n}\theta_{v i}k_{v}(\gX_{v
  i},\cdot)$ for some parameter $\gtheta \in \R^{n|\cP|}$ with
components $(\theta_{v,i}, i=1, \ldots, n, v=1, \ldots, |\cP|)$.

Let $\|\cdot\|$ denotes the usual euclidean norm in $\R^{n}$. 
For each $v \in \cP$, let $K_{v}$ be the $n \times n$ matrix with
components $(K_{v})_{i,i'} = k_{v}(\gX_{v i},\gX_{v i'})$ that satisfies $t(K^{1/2})K^{1/2}=K$.
Let $\widehat{f}_{0}$ and $\widehat{\gtheta}$ be the minimizer of the following penalized least-squares criteria:
 \begin{equation}
 C(f_{0}, \gtheta ) = \frac{1}{n}\|\gY - f_{0} 1_{n} - \sum_{v\in \cP} K_{v}
 \gtheta_{v} \|^{2}
+ \frac{1}{\sqrt{n}}\sum_{v\in \cP}\gamma_{v} \|K_{v}\gtheta_{v} \| 
+  \sum_{v\in \cP} \mu_{v} \|K^{1/2}_{v}\gtheta_{v} \|.
\label{critC.eq}
\end{equation}
Then the estimator $\widehat{f}$ defined at
Equation~\eref{fchap} satisfies
\begin{equation*}
 \widehat{f}(\gx) = \widehat{f}_{0} +
\sum_{v\in \cP} \widehat{f}_{v}(\gx_{v}) \mbox{ with } 
\widehat{f}_{v}(\gx_{v}) = \sum_{i=1}^{n}\widehat{\theta}_{v, i}k_{v}(\gX_{v, i},\gx_{v}).
\end{equation*}

Because $C(f_{0}, \gtheta)$ is a convex and separable
criteria, we propose to calculate $\widehat{\gtheta}$ using a block coordinate
descent algorithm described in the following section.

Note that the estimator
$\widehat{f}$ defined at Equation~\eref{fchap} should satisfy
$\|f_{v}\|_{\cH_v} = \|K^{1/2}_{v}\gtheta_{v} \|  \leq 1$, or
generally $\|f_{v}\|_{\cH_v} \leq r_{v}$ for some positive $r_v$, see~\eref{Fr.eq}. Usually  we
have no idea of the value of this upper-bound in practice, and we
propose to remove this contraint in the optimization
procedure. Nevertheless, if one wants to consider such an additional
constraint, the problem can be solved at the price of some additional
complication, considering a Lagrangian method for example, see
Section~\ref{algoProofs.st} for more details. 

\subsection{Algorithm}
\label{algo.st}

We will assume that
for all $v \in \cP$
the matrices $K_{v}$ are strictly definite positive. If it is not the
case, one  modifies $K_{v}$ by $K_{v} + \xi I_{n}$ where $\xi$ is a
small positive value, in order to ensure positive
definiteness.

Using a coordinate descent procedure, we minimize the  criteria
$C(f_{0}, \gtheta)$ along each group $v$ at a time. At each step of the algorithm, the criteria is
minimized as a function of the current block's parameters, while the parameters
values for the other blocks are fixed to their current values. The
procedure is repeated until convergence, considering for example that
the convergence is obtained if the norm of the difference between two
consecutive solutions is  small enough. See for exemple Boyd et
al.~\cite{Boyd2011} for optimization in such context.

For the sake of simplicity, we consider the minimization of the
following criteria:
\begin{equation}
  C'(f_{0}, \gtheta ) = \|\gY - f_{0} - \sum_{v\in \cP} K_{v}
 \gtheta_{v} \|^{2} 
+ \sum_{v\in \cP}\gamma'_{v}   \|K_{v}\gtheta_{v} \|
+ \sum_{v\in \cP} \mu'_{v}  \|K^{1/2}_{v}\gtheta_{v} \|.
\label{CritAlgo.eq}
\end{equation}
Taking
$\gamma'_{v} = \sqrt{n}\gamma_{v}$ and $\mu'_{v} = n \mu_{v}$, this is
exactly the same criteria as the one defined at Equation~\eref{critC.eq}.

Let us begin with the constant term $f_{0}$. Because the penalty
function does not depend on $f_{0}$, minimizing $C'(f_{0}, \gtheta )$
with respect to $f_{0}$ for fixed values of $\gtheta$ leads to 
\begin{equation}
 f_{0} = Y_{\cdot} - \sum_{v} \sum_{i=1}^{n} \left(K_{v}
  \gtheta_{v}\right)_{i}/n,
\label{f0.eq}
\end{equation}
where $Y_{\cdot}$ denotes the mean of
$\gY$ and $\left(K_{v} \gtheta_{v}\right)_{i}$ denotes the $i$-th
component of $K_{v} \gtheta_{v}$. In what follows, we consider a group $v$, and fix the values of the
parameters for all the other groups.  We describe the algorithm and postpone the proofs in Section~\ref{algoProofs.st}.

\medskip
\paragraph{\underline{Let us first consider the case
 where both $\mu'_{v}$ and $\gamma'_{v}$ are non zero}}

If $\partial C'_{v}$ denotes the subdifferential of
$C'(f_{0}, \gtheta )$ with respect to  $\gtheta_{v}$, we need to solve
$0 \in \partial C'_{v}$, which is equivalent to 
\begin{equation}
- 2  K_{v} \gR_{v} + 2 K^{2}_{v} \gtheta_{v} + \gamma'_{v} s_{v} +
\mu'_{v} t_{v} =0, \label{subdiff.eq}
\end{equation}
where 
\begin{equation*}
 \gR_{v} = \gY - f_{0} - \sum_{w \neq v} K_{w} \gtheta_{w}
\end{equation*}
and where $s_{v}$ and $t_{v}$  satisfy:
\begin{eqnarray*}
\mbox{if } \gtheta_{v} = 0 && \| K_{v}^{-1} s_v\| \leq 1, \mbox{ and }
\| K_{v}^{-1/2} t_v\| \leq 1 \\
\mbox{if } \gtheta_{v} \neq 0 &&  
s_{v} = \frac{K^{2}_{v}   \gtheta_{v}}{\| K_{v}\gtheta_{v}\|}, \mbox{ and }
 t_{v} = \frac{K_{v}   \gtheta_{v}}{\| K_{v}^{1/2}\gtheta_{v}\|}.
\end{eqnarray*}
The first task is to obtain necessary and sufficient conditions for
which the solution $\gtheta_{v}=0$ is the optimal one. 
Let 
\begin{equation}
J(t) = \| 2 \gR_{v}  - \mu'_{v} K_{v}^{-1} t\|^{2}, \mbox{ and }
J^{*} = \argmin \left\{ J(t), \mbox{ for } t \in \R^{n} \mbox{ such that } 
\| K_{v}^{-1/2} t\| \leq 1\right\}.
\label{RdgPb.eq}
\end{equation}
Then the  solution to
Equation~\eref{subdiff.eq} is zero if and only if $J^{*} \leq
\gamma'^{2}_{v}$. Calculating $J^{*}$  is a ridge
regression problem that can be easily solved (see
Propositions~\ref{CNS.pr} and~\ref{RdgPb.pr} in Section~\ref{algoProofs.st}).
 
If the solution to Equation~\eref{subdiff.eq} is not $\gtheta_{v} = 0$, the problem is to solve the subgradient
equation:
 \begin{equation}
\gtheta_{v} = \left( 
\frac{\mu'_{v}}{2\| K_{v}^{1/2}\gtheta_{v}\|} I_{n}
+ K_{v} + \frac{\gamma'_{v}}{2\| K_{v}\gtheta_{v}\|} K_{v}
\right)^{-1}  \gR_{v}. \label{subdiffNzero.eq}
\end{equation}
Because $\gtheta_{v}$ appears in both sides of the equation, a
numerical procedure is needed (see Proposition~\ref{Solvteta.pr}). 

\medskip
\paragraph{\underline{Other cases}}
\begin{itemize}
\item If all the $\mu'_{v}$ are equal to 0, the parameters $\gtheta_v$ are not
identifiable.
\item If all the $\gamma'_{v}$ are equal to 0, then we have to
solve a classical group-lasso problem with respect to the parameters
$\gtheta'_{v}$ defined as 
$\gtheta'_{v} = K_{v}^{1/2}\gtheta_{v}$ for all $v \in \cP$.
\item Let $v$ such that $\mu'_{v}=0$, $\gamma'_{v} \neq 0$ and assume
  that at least one of the
$\mu'_{w}$ is non zero for $w \in \cP$, $w \neq v$. Then it is shown in
Proposition~\ref{Casemu0.pr} 
that $\gtheta_{v}=0$ if and only if $2
\|R_v\| \leq \gamma'_{v}$.
\item In the same way, if $\gamma'_{v}=0$, and
$\mu'_{v} \neq 0$, then $\gtheta_{v}=0$ if and only if $2\|K_{v}^{1/2}R_v\| \leq \mu'_{v}$.
\end{itemize}

\paragraph{}
Finally the algorithm is the following:
\begin{enumerate}
\item Start with an initial value $\gtheta = \gtheta_{0}$. 
\item Calculate
  $f_{0}$ using Equation~\eref{f0.eq}. For the group $v$, calculate $\gR_{v}$ and determine if the group
  $v$ should be excluded or not either by solving the problem defined at
  Equation~\eref{RdgPb.eq} if $\mu'_{v}\neq 0$ and  $\gamma'_{v} \neq
  0$, or directly if one of them equals 0. If it is the case, set $\gtheta_{v}=0$. If
  not, solve Equation~\eref{subdiffNzero.eq} to obtain $\gtheta_{v}$. 
\item Iterate step 2. over all the groups $v$.
\item Iterate step 2. and 3. until convergence.
\end{enumerate}

\subsection{Choice of the tuning parameters}
\label{ChTP.st}
For each value of the tuning parameters $(\mu'_{v}, \gamma'_{v}), v
\in \cP$, the algorithm provides 
an estimate of $m$ and of the Sobol
indices. The problem for choosing these parameters values is crucial. We
propose to restrict this choice by considering tuning parameters
proportionnal to known weights: for all $v \in \cP$, $\mu'_{v}=\mu
\omega_{v}$ and $\gamma'_{v} = \gamma \zeta_{v}$, where the weights
$\omega_{v}$ and $\zeta_{v}$ are  fixed. For example, one can take
weights that increase with the cardinal of $v$ in order to favour
effects with small interaction order between variables. Or, according
to the theoretical result given at Theorem~\ref{oracle}, we can choose
$\omega_{v} = \widehat{\nu}^{2}_{n,v}$ and $\zeta_{v} =
\widehat{\nu}_{n,v}$, where $\widehat{\nu}_{n,v}$ is an estimate of
$\nu_{n,v}$ based on the eigenvalues of the matrix $K_{v}$. Any other
choice, depending on the problem of interest, may be relevant.

Once the weights are chosen, we estimate  $m$, on the basis of a  learning data set
$(Y_{i}, \gX_{i}), i=1, \ldots n$, for a
grid of values of $(\mu, \gamma)$. We first set
$\gamma=0$, and calculate $\mu_{\max}$ the smallest value of $\mu$ such
that the solution to the  minimization of 
\begin{equation*}
 \|\gY - f_{0} - \sum_{v\in \cP} K_{v}
 \gtheta_{v} \|^{2} 
+ \mu \sum_{v\in \cP} \omega_{v}  \|K^{1/2}_{v}\gtheta_{v} \|,
\end{equation*}
is $\gtheta_{v}=0$ for all $v \in \cP$. Then we can consider
$\mu_{\ell} = \mu_{\max} 2^{-\ell}$ for $\ell \in \left\{1, \ldots,
  \ell_{\max}\right\}$, as a grid of values for $\mu$. The grid of
values for $\gamma$ may be chosen after few attempts.

For choosing the final estimator, say $\widehat{f}$, we suppose that
we have at our disposal a testing data set $(Y^{T}_{i},
\gX^{T}_{i}), i=1, \ldots n^{T}$, and we propose two procedures. 
\begin{description}
\item[Proc. GS]
The first one uses the testing data set
 for estimating the prediction
error. Precisely, for each value of $(\mu, \gamma)$ in the grid, let
$\widehat{f}_{(\mu, \gamma)}(\cdot)$ be the estimation of $m$ obtained
with the learning data set. Then  
\begin{equation*}
\PE(\mu, \gamma) = \frac{1}{n^{T}} \sum_{i=1}^{n^{T}}
\left(Y^{T}_{i} - \widehat{f}_{(\mu, \gamma)}(\gX^{T}_{i})\right)^{2}
\end{equation*}
estimates the prediction error, and we propose to choose the pair $(\mu, \gamma)$ that
minimizes $\PE(\mu, \gamma)$, say $(\widehat{\mu},
\widehat{\gamma})$. 
Finally the estimator, denoted $\widehat{f}^{\T}$ is defined as
$\widehat{f}^{\T} = \widehat{f}_{(\widehat{\mu},
  \widehat{\gamma})}$. In the following, we will refer to this
procedure as the Group-Sparse procedure.

\item[Proc. rdg]
Doing the parallel with the
inconsistency of the lasso for estimating the support of the
parameters in the classical
regression problem, we propose to choose the tuning parameter that minimizes the risk of the ridge estimator over the support estimated
by  the
ridge-group-sparse procedure.  Indeed, if the tuning parameter  is chosen to minimize the prediction error, the lasso is not consistent for support estimation (see~\cite{LengLiWahba} for example). One idea to overcome this
problem, is to choose the tuning 
parameter that minimizes the risk of the 
Gauss-lasso estimator which is calculated in two steps: For a given value of the tuning parameter, the
estimation of the support of the parameter is estimated using a lasso
procedure, then the least square estimator over this support is
calculated. When the objective is support estimation, some
numerical simulations~\cite{RigolletTsybakov2011} and theoretical results ~\cite{JavanmardMontanari2013} suggest that it may be more advisable not to apply the selection schemes based on
prediction risk  to the
lasso estimators, but rather to the Gauss-lasso estimators.
Our procedure, called the ridge procedure,  applies the same idea  in the framework
of sparse nonparametric estimation. Precisely
it  considers the collection of supports composed of
the different $S_{\widehat{f}_{(\mu, \gamma)}}$ when $(\mu, \gamma)$
belongs to the grid. For each support $S$ in this collection, we
estimate $f$ using a ridge procedure assuming that the support of $f$
is $S$: for a given $\lambda$,
$f^{\rdg}_{\lambda, S}$ is defined as follows:
\begin{equation*}
 f^{\rdg}_{\lambda, S} = \argmin\left\{ 
\frac{1}{n} \sum_{i=1}^{n} \left(
Y_{i} - f_{0} - \sum_{v \in S} f_{v}(\gX_{v, i}) \right)^{2}
+ \lambda \sum_{v \in S} \|f_{v}\|^{2}_{\cH_{v}}, f=f_{0} + \sum_{v\in
S} f_{v}, f_{v} \in \cH_{v}
\right\}.
\end{equation*}
We choose a grid of values for $\lambda$, and for each $S$ in the
collection and $\lambda$ in the grid, we estimate the prediction error
$\PE(\lambda, S)$ defined as follows:
\begin{equation*}
 \PE(\lambda, S) = \frac{1}{n^{T}} \sum_{i=1}^{n^{T}}
\left(Y^{T}_{i} - {f}^{\rdg}_{\lambda, S}(\gX^{T}_{i})\right)^{2}.
\end{equation*}
For each $S$ in the collection, let $\widehat{\lambda}(S)$ be the
minimizer of $\PE(\lambda, S)$ when $\lambda$  varies in the grid, and
let $\widehat{S}$ be the minimizer of $\PE(\widehat{\lambda}(S), S)$,
then the estimator denoted $\widehat{f}^{\rdg}$ is defined as $\widehat{f}^{\rdg} =
f^{\rdg}_{\widehat{\lambda}(\widehat{S}),\widehat{S}}$.

\end{description}
If a testing data set is not available, we can use the classical
  V-fold cross validation (see~\cite{ArlotCelisse2010} for example) either to estimate $\PE(\mu, \gamma)$ or to
  estimate  $\PE(\lambda, S)$.

\section{Simulation study}
\label{simul.st}

In order to evaluate the performances of our method for estimating a
meta-model and 
sensitivity indices of a function $m$ we carried out a simulation
study. We consider the $g$-function of Sobol defined on $[0,1]^{d}$
as
\begin{equation*}
m(\gx) = \prod_{a=1}^{d} \frac{|4 x_{a} -2| + c_{a}}{1 + c_{a}}, c_{a} >0,
\end{equation*}
whose Sobol indices can be expressed analytically (see Saltelli et
al.~\cite{Saltelli2000}). Following the simulation experiment proposed
by Durrande et 
al.~\cite{Durrande2013}, we take $d=5$ and $(c_1, c_2, c_3, c_4,
c_5) = (0.2, 0.6, 0.8, 100, 100)$. The lower the value of $c_{a}$, the
more significant the variable $x_{a}$. 
The variables $X_{a}, a=1, \ldots,
d$ are independent and uniformly distributed on $[0,1]$. We consider
the regression model $Y_{i}  = m(\gX_{i}) + \sigma \varepsilon_{i}$,
for $i=1, \ldots, n$, with $\cN(0,1)$ independent error terms
$\varepsilon_{i}$. 

\paragraph{\underline{Simulation design}}

We present the results for $n \in \left\{50, 100, 200\right\}$,
$\sigma \in \left\{0,  0.2 \right\}$. For all $a=1, \ldots, d$, the kernels
$k_{a}$ are the same: we considered the
Brownian kernel, $k^{b}(x, x') = 1 + \min\left\{x,x'\right\}$, the
Mat\'ern kernel, $k^{m}(x, x') = (1+2|x-x'|)\exp(-2|x-x'|)$, and the
Gaussian kernel, $k^{g}(x, x') = \exp(x-x')^{2}$. 

For each simulation, we generate three independent data sets as follows
: a Latin Hypercube Sample 
of the inputs is simulated to give the matrix $\gX$ with $n$ rows and
$d=5$ columns, and a $n$-sample of independent centered Gaussian
variable with variance 1 is simulated. This operation is repeated
three times  in order to obtain the learning and testing data sets
and a third data set for estimating the estimators performances. As
explained in Section~\ref{ChTP.st}, we
choose optimal values of the tuning parameters, $(\mu,
\gamma)$ by minimizing a prediction error $\PE$, and get an
estimator $\widehat{f}$ of
$m$, as well as estimates of the Sobol indices:
\begin{equation*}
\widehat{S}_{v} = 
\frac{\widehat{\Var}\left(m_{v}(\gx_{v})\right)}{
\sum_{w \in \cP}\widehat{\Var}\left(m_{w}(\gx_{w})\right)}.
\end{equation*}
Let $\Omega_{v}$ be the matrix whose components 
satisfy 
\begin{equation*}
 (\Omega_{v})_{i,i'} = \prod_{a\in v} 
\E_{U \sim P_{a}} \left(k_{0a}(U,X_{a,i}) k_{0a}(U, X_{a, i'}) \right)
\mbox{ for } i,i' =
1, \ldots, n.
\end{equation*}
The estimator $\widehat{\Var}\left(m_{v}(\gx_{v})\right)$ is
calculated as follows: 
\begin{equation*}
 \widehat{\Var}\left(m_{v}(\gx_{v})\right) 
 =  \widehat{\alpha}_{v}^{T} \Omega_{v} \widehat{\alpha}_{v}.
\end{equation*}

We also propose to estimate these quantities by their empirical
variances as in Equation~\eref{Vemp.eq}.

\medskip
\paragraph{\underline{Performance indicators}}
To evaluate the performances of our method for estimating a
meta-model, we use the
classical coefficient of determination $R^{2}$ estimated using the
third data set $(Y_{i}^{P}, \gX^{P}_{i}), i=1, \ldots, n $:
\begin{equation*}
 R^{2} = 1 - \frac{\sum_{i=1}^{n}
\left(Y^{P}_{i} - \widehat{f}(\gX^{P}_{i})\right)^{2}}{
\sum_{i=1}^{n}\left(Y^{P}_{i} - Y^{P}_{\cdot}\right)^{2}}.
\end{equation*}
Moreover we calculate the empirical risk $\ER = \| m - \widehat{f}\|_{n}^{2}$.
 For each simulation $s$, we get $R^{2}_{s}$ and $\ER_{s}$ and we report the means  of
these quantities over all simulations.

Similarly, for each $v$, and each simulation $s$, we get
$\widehat{S}_{v, s}$
and we report its mean, $ \widehat{S}_{v,
  \cdot}$, its estimated standard-error, and to sum up the behaviour of our
procedure for estimating the sensitivity indices, we estimate the 
global error, denoted  $\GE$, defined as follows
\begin{eqnarray*}
 \GE & = & \sum_{v} (\widehat{S}_{v, s} -
 S_{v})^{2}.
\end{eqnarray*}

In order to assess the performances of our procedure for selecting the
greatest Sobol Indices, precisely those that are greater than some
small quantity as $\rho=10^{-4}$,  we
calculate for each group $v \in \cP$, the percentage of simulations
for which  $v$ is in the support of the estimator $\widehat{f}$. Then
we average these quantities, on one hand over groups $v$ such that
$S_{v} > \rho$, and on the other hand,  over groups $v$ such that
$S_{v} \leq \rho$. Let us
denote these quantities ${\mathrm{pSel}}_{S_{v} > \rho}$ and ${\mathrm{pSel}}_{S_{v} \leq \rho}$ respectively.

\medskip
For each $(\mu,
\gamma)$ in a grid of values, the estimator
$\widehat{f}_{\mu, \gamma}$ is
defined as the minimizer of the criteria given at
Equation~\eref{CritAlgo.eq} taking 
$\omega_{v}=\zeta_{v}=1$ for all $v \in \cP$. In order to save
computation time,  we restrict the optimisation to sets $v$ such that $|v| \leq
3$. Some preliminary simulations showed that the terms
coresponding to $|v| \geq 4$ are nearly always equal to 0.

\medskip
\paragraph{\underline{Choosing the tuning parameters}}

Let us begin with the comparison of the two methods proposed for
choosing the final estimator, see Section~\ref{ChTP.st}. The results
are given in Tables~\ref{R2.tb} and~\ref{ER.tb}. It appears that the procedure based
on the ridge estimator after selection of the groups outperforms the
method based on the group-sparse estimator. As expected both methods 
perform better when $n$ increases, and when $\sigma=0$.

\begin{table}
\begin{center}
\begin{tabular}{l|ccc|ccc}
& \multicolumn{3}{c|}{$\sigma=0$}& \multicolumn{3}{c}{$\sigma=0.2$} \\
&                         $n=50$ &$n=100$&$n=200$& $n=50$&$n=100$&$n=200$\\
Proc. GS       &0.814&0.920 &0.959&0.737&0.835&0.889\\ 
Proc. rdg      &0.874&0.976 &0.989&0.763&0.854&0.892\\
\end{tabular}
\end{center}  

\caption{\label{R2.tb}  Estimated coefficient of determination $R^{2}$
  for different values of $n$ 
  and $\sigma$, with the Mat\'ern kernel.}
\end{table}

\begin{table}
\begin{center}
\begin{tabular}{l|ccc|ccc}
& \multicolumn{3}{c|}{$\sigma=0$}& \multicolumn{3}{c}{$\sigma=0.2$} \\
&                         $n=50$ &$n=100$&$n=200$& $n=50$&$n=100$&$n=200$\\
Proc. GS   & 0.033& 0.0137 & 0.0139& 0.051& 0.028& 0.020\\ 
Proc. rdg  & 0.011 & 0.0009& 0.0007& 0.042& 0.022 & 0.013\\
\end{tabular}
\end{center}  
\caption{\label{ER.tb} Estimated empirical risk $\ER$ for different values of $n$
  and $\sigma$ with the Mat\'ern kernel. }
\end{table}

Similarly, the Sobol indices are better estimated, in the sense of the
global error, with the procedure based on the ridge estimate of the
metamodel, see Table~\ref{GE.tb}. The means of the estimators for the
Sobol indices greater than $\rho=10^{-4}$ are given in
Tables~\ref{SI.n50.s0.tb} and~\ref{SI.n100.s02.tb}. It appears than
$S_{\{1\}}$ is over-estimated using the procedure based on the
group-sparse estimator, leading to under-estimate the Sobol indices
associated with interactions of order 2. This tendancy is much less
pronounced with the procedure based on the
ridge estimator.

\begin{table}
\begin{center}
\begin{tabular}{l|ccc|ccc}
& \multicolumn{3}{c|}{$\sigma=0$}& \multicolumn{3}{c}{$\sigma=0.2$} \\
&                         $n=50$&$n=100$&$n=200$& $n=50$&$n=100$&$n=200$\\
Proc. GS    & 1.91 & 0.79  & 0.45  & 2.41  & 1.16  & 0.54\\ 
Proc. rdg   & 0.80 & 0.10  & 0.03  & 1.50  & 0.47  & 0.15\\
\end{tabular}
\end{center}
\caption{\label{GE.tb}  Estimated global error $\GE \times 100$ for different values of $n$
  and $\sigma$ with the Mat\'ern kernel. }
\end{table}

\begin{table}
\begin{center}
\begin{tabular}{l|ccccccc|c}
  &$v=\{1\}$&$v=\{2\}$&$v=\{3\}$&$v=\{1,2\}$&$v=\{1,3\}$&$v=\{2,3\}$&$v=\{1,2,3\}$&
  sum \\ \hline
 S.I.                    &43.3      &24.3      &19.2      &5.63       &4.45       &2.50      &0.579        &99.98\\
Proc. GS      &50.1 (6.2)&26.5 (5.4)&20.9 (4.9)&0.69 (0.8) & 0.63 (1.1)&0.51 (0.9)&0.02 (0.07) &99.29\\
Proc. rdg     &45.4 (4.3)&25.3 (3.4)&20.3 (3.5)&3.08 (2.5)  &2.18 (2.1)&1.44 (1.8)&0.09 (0.5)&98.66 \\
\end{tabular}
\end{center}
\caption{\label{SI.n50.s0.tb} The first line of the table gives the true values
  of the Sobol indices $\times 100$ greater than $10^{-2}$, as well as  their sum in
  the last columns. The following lines give
  the mean of the estimators as well as their standard-error (in
  parenthesis), calculated over 100 simulations, for $n=50$,
  and $\sigma=0$ with the Mat\'ern kernel.}
\end{table}

\begin{table}
\begin{center}
\begin{tabular}{l|ccccccc|c}
  &$v=\{1\}$&$v=\{2\}$&$v=\{3\}$&$v=\{1,2\}$&$v=\{1,3\}$&$v=\{2,3\}$&$v=\{1,2,3\}$&
  sum \\ \hline
 S.I.                    &43.3      &24.3      &19.2      &5.63      &4.45      &2.50      &0.579      &99.98\\
Proc. GS       &47.5 (5.4)&26.2 (4.9)&19.8 (3.7)&2.35 (1.5)&1.45 (1.2)&0.84 (0.8)&0.03 (0.1)&98.95\\
Proc. rdg    &43.5 (3.9)&25.0 (3.6)&19.7 (2.8)&4.85 (1.7)&3.39 (1.5)&2.02 (1.3)&0.05 (0.3)&99.02\\
\end{tabular}
\end{center}
\caption{\label{SI.n100.s02.tb}  The first line of the table gives the true values
  of the Sobol indices $\times 100$ greater than $10^{-2}$, as well as  their sum in
  the last columns. The following lines give
  the mean of the estimators as well as their standard-error (in
  parenthesis), calculated over 100 simulations, for $n=100$,
  and $\sigma=0.2$ with the Mat\'ern kernel.}
\end{table}

Let us now consider the performances of the procedure for selecting
the non zero Sobol indices. In Tables~\ref{pSel.n50.s0.tb}
and~\ref{pSel.n100.s02.tb} we report the percentages of simulations
for which the Sobol indices smaller (respectively greater) than $\rho$ are
selected, and for which each of Sobol index greater than $\rho$ is
selected. From these results, we conclude that the procedure based on
the ridge estimator is more strict for selecting non-zero Sobol
indices.

\begin{table}
\begin{center}
\begin{tabular}{l|ccccccccc}
&$\SI < \rho$ & $\SI\geq \rho$   &$v=\{1\}$&$v=\{2\}$&$v=\{3\}$&$v=\{1,2\}$&$v=\{1,3\}$&$v=\{2,3\}$&$v=\{1,2,3\}$
\\ \hline
Proc. GS       & 17.9& 68& 100&100&100&72&66&67&9 \\
Proc. rdg      & 6.3 & 51& 100&100&100&72&62&47&4 \\
\end{tabular}
\end{center}
\caption{\label{pSel.n50.s0.tb}The first two columns
  give respectively ${\mathrm
  pSel}_{S_{v} \leq \rho}$ and ${\mathrm pSel}_{S_{v} >
  \rho}$. The last columns give the values of ${\mathrm
  pSel}_{v}$ for each group $v$ such that $S_{v} > \rho$. Results for
$n=50$ and $\sigma=0$ with the Mat\'ern kernel. }
\end{table}

\begin{table}
\begin{center}
\begin{tabular}{l|ccccccccc}
&$\SI < \rho$ & $\SI\geq \rho$   &$v=\{1\}$&$v=\{2\}$&$v=\{3\}$&$v=\{1,2\}$&$v=\{1,3\}$&$v=\{2,3\}$&$v=\{1,2,3\}$
\\ \hline
Proc. GS       &38 & 85& 100&100&100&100& 99& 98& 18 \\
Proc. rdg      & 6 & 58& 100&100&100& 98& 92& 77&  3 \\
\end{tabular}
\end{center}
\caption{\label{pSel.n100.s02.tb} blablaThe first two columns
  give respectively ${\mathrm
  pSel}_{S_{v} \leq \rho}$ and ${\mathrm pSel}_{S_{v} >
  \rho}$. The last columns give the values of ${\mathrm
  pSel}_{v}$ for each group $v$ such that $S_{v} > \rho$. Results for
$n=100$ and $\sigma=0.2$ with the Mat\'ern kernel.}
\end{table}

\medskip
\paragraph{\underline{Comparing different kernels}}

Finally we compare the performances of the procedures for different
kernels, see Tables~\ref{Ckernel.ER.tb} and~\ref{Ckernel.GE.tb}. The
means of the estimated empirical risk, $\ER$, and of the global error
for estimating the sensitivity indices, $\GE$, are calculated for each
kernel.  It appears that the Mat\'ern kernel gives the best results,  except for the case
 $n=50$ and $\sigma=0$ where the empirical risk of
 $\widehat{f}^{\T}$ is smaller for the Brownian kernel.  

In practice,
 one may want to choose the kernel according to the smallest
 prediction error. For that purpose, we propose  to calculate the estimators
 $\widehat{f}^{\T}$  and/or $\widehat{f}^{\rdg}$ for each kernel, as described in
 Section~\ref{ChTP.st}, as well as  their associated prediction
 errors. Then we choose 
 the kernel for which the prediction error is minimized. The results
 are reported under the column ``mixed'' in
 Tables~\ref{Ckernel.ER.tb} and~\ref{Ckernel.GE.tb}. It appears that
 the estimated empirical risks for this ``mixed'' procedure are nearly
 equal to the minimum estimated empirical risks over the different kernels.

\begin{table}
\begin{center}
\begin{tabular}{l|cccc|cccc}
& \multicolumn{4}{c|}{$n=50$, $\sigma=0$}&\multicolumn{4}{c|}{$n=100$, $\sigma=0.2$} \\ \hline
               &Mat\'ern&Gaussian&Brownian&mixed&Mat\'ern&Gaussian&Brownian&mixed \\\hline
Proc. GS       & 0.033& 0.054  & 0.027  &0.028&0.033 &0.054   &0.027   &0.028  \\
Proc. rdg      & 0.011& 0.024  & 0.025  &0.011&0.011 &0.024   &0.025   &0.023\\
\end{tabular}
\caption{\label{Ckernel.ER.tb} Estimated empirical risk $\ER$: Performances of the procedures according to the kernel choice.}
\end{center}
\end{table}

\begin{table}
\begin{center}
\begin{tabular}{l|cccc|cccc}
& \multicolumn{4}{c|}{$n=50$, $\sigma=0$}&\multicolumn{4}{c|}{$n=100$, $\sigma=0.2$} \\ \hline
            &Mat\'ern&Gaussian&Brownian&mixed&Mat\'ern&Gaussian&Brownian&mixed \\\hline
Proc. GS    &1.91 &2.49    &2.19     &1.88 &1.16  &1.90    &1.69&1.16 \\
Proc. rdg   &0.80 &1.39    &1.40     &0.85 &0.48  &0.83    &0.73&0.51\\
\end{tabular}
\caption{\label{Ckernel.GE.tb}  Estimated global error $\GE \times
  100$: Performances of the procedures according to the kernel choice.}
\end{center}
\end{table}

\section{\label{sketch.st}Sketch of proof of Theorem~\ref{oracle}}
We give here a
sketch of the proof and we postpone to 
Section~\ref{proofs.st} for complete statements. In particular, we  denote by $C$ constants that
vary  from an equation to the other, and we assume that $\sigma=1$.

The proof of Theorem~\ref{oracle} starts in the same way as the proof
of Theorem 1 in  Raskutti et al.~\cite{RWY}. Nevertheless
it differs in several points, in particular because the terms occuring in the decomposition of functions in $\cH$ depend
on several variables and thus are not independent. Indeed, $f_v(\gX_{v})$ and $f_v'(\gX_{v'})$ are not independent
  as soon as the groups  $v$ and $v'$ share some of the variables
  $X_{a}, a=1, \ldots, d$. Moreover, we do not
assume that the function $m$ is in $\cF$.

Starting from the definition of $\widehat{f}$, some simple calculation
(see Equation~\eref{base2.eq})
give that for all $f\in \cF$ 
\begin{equation*}
 C \|m - \widehat{f}\|_{n}^{2} \leq \|m - f\|_{n}^{2} +
\left| \frac{1}{n}\sum_{i=1}^{n} \varepsilon_{i} (\widehat{f}(\gX_{i})-f (\gX_{i}))\right|
+ \sum_{v \in
  S_{f}}\left( \gamma_{v} \|\widehat{f}_{v}-f_{v}\|_{n} + \mu_{v}
\|\widehat{f}_{v} - f_{v}\|_{\cH_{v}}\right).
\end{equation*}
If we set $g=\widehat{f} - f$, then $g \in
\cH$, $g=g_{0} + \sum_{v} g_{v}$, with $g_{v} = \widehat{f}_v - f_v$, and for each $v$, $\|g_{v}\|_{\cH_{v}}
\leq 2$.

The main problem is now to control the empirical process.  
For each $v$, letting $\lambda_{n,v}$ as in (\ref{lambda.eq}),
we state (see Lemma~\ref{Tau}, page~\pageref{Tau}) that, with high
probability, 
\begin{eqnarray*}
\mbox{if }  \|g_{v}\|_{n} \leq \lambda_{n,v} \|g_{v}\|_{\cH_{v}} &\mbox{
  then } &
\left| \sum_{i=1}^{n}\varepsilon_{i} g_{v}(\gX_{v, i})\right| \leq C n \lambda_{n,v}^{2}
\|g_{v}\|_{\cH_{v}} \\
\mbox{if }  \|g_{v}\|_{n} > \lambda_{n,v} \|g_{v}\|_{\cH_{v}} &\mbox{
  then } &
\left| \sum_{i=1}^{n}\varepsilon_{i} g_{v}(\gX_{v, i})\right| \leq C n \lambda_{n,v}
\|g_{v}\|_{n}.
\end{eqnarray*}

Therefore, if for all $v$, $\mu_{v}$ and  $\gamma_{v}$ satisfy
Equation~\eref{condmu}, we deduce that with high probability (setting $g=\widehat{f}-f$)
\begin{equation*}
C \|m - \widehat{f}\|_{n}^{2} \leq   \|m - f\|_{n}^{2}
\sum_{v \in S_f} \left( \gamma_{v}\|g_{v}\|_{n} + \mu_{v}
\|g_{v}\|_{\cH_{v}} \right)
+ \sum_{v \notin S_{f}}
\left(\gamma_{v} \|\widehat{f}_{v}\|_{n} + \mu_{v} \|\widehat{f}_{v}\|_{\cH_{v}}\right).
\end{equation*}
Besides we can express the decomposability property of the penalty as
follows (see lemma~\ref{lemme2}, page~\pageref{lemme2}): with high probability (in the set
where the empirical process  is controled as stated above), 
\begin{equation*}
\sum_{v \notin
  S_{f}}\left(\gamma_{v} \|\widehat{f}_{v}\|_{n} + \mu_{v}
\|\widehat{f}_{v}\|_{\cH_{v}}\right)
\leq  C \sum_{v \in
  S_{f}}\left(\gamma_{v} \|g_{v}\|_{n} + \mu_{v}
\|g_{v}\|_{\cH_{v}}\right).
\end{equation*}
Putting the things together, and noting again that $\|g_{v}\|_{\cH_{v}} \leq
2$, we obtain the following upper bound
\begin{equation*}
 C \|m - \widehat{f}\|_{n}^{2} \leq  \|m - f\|_{n}^{2} + \sum_{v \in S_{f}} \left(\ \mu_{v} +
\gamma_{v}  \|g_{v}\|_{n}\right).
\end{equation*}
The last important step consists in comparing $\sum_{v \in S_{f}} \|g_{v}\|_{n} $ to
$\|\sum_{v \in S_{f}} g_{v}\|_{n}$. More precisely, it can be shown
(see lemma~\ref{norm2normn}, page~\pageref{norm2normn}) that for all
$v\in \cP$, with high probability, we have
\begin{equation*}
  \|g_{v}\|_{n}  \leq 2 \|g_{v}\|_{\bL^{2}(\PX)} + \gamma_{v}.
\end{equation*}
Using the orthogonality assumption between the spaces $\cH_{v}$, we
have  $\sum_{v \in S_{f}}
\|g_{v}\|^{2}_{\bL^{2}(\PX)} = \|\sum_{v}
g_{v}\|^{2}_{\bL^{2}(\PX)}$, and thus we get
\begin{equation*}
 C \|m - \widehat{f}\|_{n}^{2} \leq  \|m - f\|_{n}^{2}+ \sum_{v \in S_{f}} \mu_{v}
  + \sum_{v \in S_{f}} \gamma^{2}_{v} + 
\| \sum_{v \in    S_{f}} 
(\widehat{f_{v}} - f_{v})^{2}\|^{2}_{\bL^{2}(\PX)}.
\end{equation*}
Finally it remains to consider different cases according to the
rankings of $\|\widehat{f}-f\|^{2}_{\bL^{2}(\PX)}$,
$\|\widehat{f}-f\|^{2}_n$ and $\sum_{v \in
    S_{f}} \mu_{v} + \gamma^{2}_{v}$ to get the result of Theorem~\ref{oracle}.

\section{Proofs}
\label{proofs.st}

Recall that we cconsider the regression  model defined at
Equation~\eref{RegMod.eq}, where $\gX$ has distribution $\PX =
P_{1}\times\ldots\times P_{d}$ defined on $\cX$ a compact subset of $\R^{d}$, and
$\varepsilon$ is  distributed as $\cN(0, \sigma^{2})$. We denote
by $\PXe$ the distribution of $(\gX, \geps)$. We observe a $n$ sample $(Y_i,
\gX_i), i=1, \ldots, n$ with law $\PXe$. 

The notation and the procedure are given in
Sections~\ref{procedure.st} and~\ref{ThRes.st}.

Let us add on few notations that will be used along the
proofs. 

For $v \in \mathcal{P}$ we denote $|v|$ the cardinal of $v$. 
For a function $\phi :\mathbb{R}^{|v|}\mapsto \mathbb{R}$,  we denote
$V_{n,\varepsilon}$ the empirical process defined by
\begin{eqnarray}
\label{procemp}
V_{n,\varepsilon}(\phi)=\frac{1}{n}\sum_{i=1}^n \varepsilon_i \phi(\gX_{v,i}).
\end{eqnarray}
For the sake of simplicity we assume  $\sigma=1$. Moreover, we set
$R'=1$,
see~\eref{supfv.eq}. Consequently,  for any
function $f \in \cH$, $\|f_v\|_{\cH_v} \leq 1$, and $\|f_v\|_{\infty}
\leq \|f\|_{\cH_{v}}$. The proofs can be done exactly in the same way by
considering the general case.
In the proofs, the $\|\cdot \|_{\bL^{2}(\PX)}$ norm will be denoted by $\| \cdot \|_{2}$.

\subsection{Proof of  Theorem \ref{oracle}}
The proof  is based on four main  lemmas proved in
Section~\ref{proofsintermediate}. In Section~\ref{intermediate} other
lemmas used all along the proof are stated. Their proof are postponed
to Section~\ref{ProofsIntLemm.st}.

Let us first establish inequalities that will be used in the following.
Let $f \in \cH$ and  $v\in S_f$ (see~\eref{sparsity}). 

Using that for any $v \in S_f$, and any norm $\|\cdot
\|$ in $\cH_v$, $\| f_v\|-\| \widehat{f}_v\|
\leq \| f_v-\widehat{f}_v\|$ and that for any $v \notin
S_f$, $\| f_v\|=0$, 
we get that
\begin{eqnarray}
\sum_{v \in \mathcal{P}} \mu_v\| f_v\|_{\mathcal{H}_v}-\sum_{v \in \mathcal{P}} \mu_v\| \widehat{ f}_v\|_{\mathcal{H}_v}
\leq \sum_{v\in S_f} \mu_v\| f_v-\widehat{f}_v\|_{\mathcal{H}_v}-\sum_{v\in S_f^c}
\mu_v\| \widehat{
  f}_v\|_{\mathcal{H}_v},\label{decomppen_H} \\
\sum_{v \in \mathcal{P}} \gamma_v\| f_v\|_{n}-\sum_{v \in \mathcal{P}}  \gamma_v\| \widehat{ f}_v\|_{n}
\leq \sum_{v\in S_f} \gamma_v\| f_v-\widehat{f}_v\|_{n}-\sum_{v\in S_f^c}
\gamma_v\| \widehat{ f}_v\|_{n}. \label{decomppen_n}
\end{eqnarray}

Combining  \eref{decomppen_H},  and 
\eref{decomppen_n},  to the fact that for any function  $f\in \mathcal{H}$, $\mathcal{L}(\widehat{f})\leq \mathcal{L}(f)$, we obtain that
\begin{equation*}
 \| m-\widehat{f}\|_n^2 \leq \| m-f\|_n^2+ B
\end{equation*}
with
\begin{equation}
B=  2V_{n,\varepsilon}\big(\widehat{f}-f\big)+
 \sum_{v\in S_f} \left[\mu_v\| \widehat{f}_v-f_v\|_{\mathcal{H}_v}+\gamma_v\| \widehat{f}_v-f_v\|_{n}\right]
- \sum_{v\in S_f^c}
\left[\mu_v\| \widehat{ f}_v\|_{\mathcal{H}_v}+  \gamma_v\| \widehat{f}_v\|_{n}\right].\label{B.eq}
\end{equation}
If $\| m-f\|_n^2 \geq B$,  we immediately get the result
since in that case
\begin{equation*}
 \| m-\widehat{f}\|_n^2 \leq 
2 \| m-f\|_n^2 \leq 
2 \| m-f\|_n^2 + \sum_{v\in  S_f}\mu_v+\sum_{v\in S_f}\gamma^2_v.
\end{equation*}
If  $\| m-f\|_n^2 < B$, we get that
\begin{eqnarray}
\|
\widehat{f}-m\|_n^2 &\leq & 2 B \label{base}\\
&\leq & 4 \vert V_{n,\varepsilon}\big(\widehat{f}-f\big)\vert
+2\sum_{v\in S_f}\left[ \mu_v \|
  \widehat{f}_v-f_v       \|_{\mathcal{H}_v}+\gamma_v \|
  \widehat{f}_v-f_v       \|_{n}\right].
\label{base2.eq}
\end{eqnarray}

The control of the empirical process $\vert
V_{n,\varepsilon}\big(\widehat{f}-f\big)\vert$ is given by the
following lemma (proved in Section~\ref{ProofTau}, page~\pageref{ProofTau}).
 \begin{lem} 
\label{Tau}
Let $V_{n,\varepsilon}$ be defined in \eref{procemp}.
For any $f$ in $\cF$, we
consider the event   $\mathcal{T}$ defined as
\begin{eqnarray}
\label{evtTau}\mathcal{T}=\left\lbrace \forall f \in \cF, \forall
  v \in \mathcal{P} , \vert
  V_{n,\varepsilon}\big(\widehat{f}_v-f_v\big) 
\vert \leq  
\kappa \lambda_{n,v}^2\| \widehat{f}_v-f_v\|_{\mathcal{H}_v}
+ 
\kappa \lambda_{n,v}\| \widehat{f}_v-f_v\|_{n}\right\rbrace,
\end{eqnarray}
where the quantities $\lambda_{n,v}$ are defined by
Equation~\eref{lambda.eq} and where $\kappa= 10 + 4 \Delta$. Then, for some positive constants $c_1, c_2$,
$$\PXe \left(\mathcal{T}\right)\geq 1-c_1\sum_{v \in \mathcal{P}}\exp(-nc_2 \lambda_{n,v}^2).$$
\end{lem}

Conditionning on  $\mathcal{T}$, Inequality~\eref{base2.eq} becomes
\begin{eqnarray*}
 \| \widehat{f}-m\|_n^2 &\leq &
4 \kappa \sum_{v \in \mathcal{P}} \left[ \lambda_{n,v}^2 \|
\widehat{f}_v-f_v\|_{\mathcal{H}_v} 
+ \lambda_{n,v}\| \widehat{f}_v-f_v\|_{n}\right] + 2\sum_{v\in S_f} \left[\mu_v\|
  \widehat{f}_v-f_v\|_{\mathcal{H}_v}+ \gamma_v\| \widehat{f}_v-f_v\|_{n}\right],
\end{eqnarray*}
which may be decomposed as follows
\begin{eqnarray*}
\| \widehat{f}-m\|_n^2 &\leq&   
\sum_{v\in S_f} \left[
4\kappa \lambda_{n,v}^2+2\mu_v\right] \|
\widehat{f}_v-f_v\|_{\mathcal{H}_v} + \sum_{v\in
  S_f}\left[
4\kappa \lambda_{n,v}+2\gamma_v\right]\|
\widehat{f}_v-f_v\|_{n} +
\\&& 4\sum_{v\in S_f^c} 
\kappa \lambda_{n,v}^2 \| \widehat{f}_v-f_v\|_{\mathcal{H}_v} +
4\sum_{v\in S_f^c} 
\kappa \lambda_{n,v}\| \widehat{f}_v-f_v\|_{n}.
\end{eqnarray*}

If we choose $C_{1}\geq \kappa$ in Theorem~\ref{oracle}, then $\kappa \lambda_{n,v}^{2} \leq \mu_{v}/2$
and $\kappa \lambda_{n,v}\leq  \gamma_{v}/2$ and 
the previous inequality becomes
\begin{equation}
  \|  \widehat{f}-m\|_n^2  \leq 
6 \sum_{v\in S_f}\left[\mu_v\| \widehat{f}_v-f_v\|_{\mathcal{H}_v}+\gamma_v
\| \widehat{f}_v-f_v\|_{n}\right]
+
4\sum_{v\in S_f^c}\left[\mu_v \| \widehat{f}_v\|_{\mathcal{H}_v}+\gamma_v
\| \widehat{f}_v\|_{n}\right]. \label{base3}
\end{equation}

Next we use the decomposability property of the penalty expressed in the
following lemma (proved in Section~\ref{Prooflemme2} page~\pageref{Prooflemme2}).
\begin{lem}
\label{lemme2}
For any $f \in \cF$, under the assumptions of Theorem \ref{oracle}
with $C_{1} \geq \kappa$,
conditionnally on $\cT$, see~\eref{evtTau}, we have 
\begin{equation*}
 \sum_{v\in S_f^c} \mu_v \|\widehat{
   f}_v\|_{\mathcal{H}_v} +\sum_{v\in S_{f}^c}
 \gamma_v\| \widehat{ f}_v\|_{n} \leq  
3 \sum_{v\in S_f} \mu_v\|
\widehat{f}_v-f_v\|_{\mathcal{H}_v}+
3 \sum_{v\in S_f} \gamma_v\| \widehat{f}_v-f_v\|_{n}.
\end{equation*}
\end{lem}
Hence, by combining \eref{base3}  and Lemma~\eref{lemme2} we obtain 
\begin{eqnarray*}
 \|  \widehat{f}-m\|_n^2
&\leq&  18 
\sum_{v\in S_f}\big[\mu_v \| \widehat{f}_v-f_v\|_{\mathcal{H}_v} +\gamma_v \| \widehat{f}_v-f_v
\|_{n
}\big].\end{eqnarray*}
For each $v$, $\|
\widehat{f}_v-f_v\|_{\mathcal{H}_v} \leq 2$ (because the functions  $\widehat{f}_v$ et $f_v$ belong to the
class $\cF$, see~\eref{calF.eq}), and 
consequently, for some constant $C$,
\begin{equation}
 \label{borne}
 \| \widehat{f}-m\|_n^2  \leq C\left\{
 \sum_{v\in  S_f}\mu_v+
 \sum_{v\in S_f}\gamma_v \| \widehat{f}_v-f_v
\|_{n}\right\}.
\end{equation}

To finish the proof it remains to compare the two quantities $\sum_{v\in S_f} \| \widehat{f}_v-f_v
\|_{n}^2$ and  $\| \sum_{v\in S_f} \widehat{f}_v-f_v
\|_{n}^2$. For that purpose we show that 
$\| \sum_{v\in S_f}\widehat{f}_v-f_v\|_n$ is less than $ \|
\sum_{v\in S_f}\widehat{f}_v-f_v\|_2^2$ plus an
additive term coming from concentration results (see the Lemma given
below). Next, thanks to the orthogonality of the spaces $\cH_{v}$
with respect to $\bL^{2}(\PX)$, 
$ \|\sum_{v\in  S_f}\widehat{f}_v-f_v\|_2^2=\sum_{v \in  S_f} 
\|\widehat{f}_{v}-f_{v}\|_2^2$. To
conclude,  it remains to
consider several cases, according to the rankings of 
$\|\sum_{v\in  S_f}\widehat{f}_v-f_v\|_2^2$,
$\|\sum_{v\in  S_f}\widehat{f}_v-f_v\|_n^2$, and $d^{2}(f)$.
This is the subject
of the following lemma whose proof is given in
Section~\ref{Proofnorm2normn}, page~\pageref{Proofnorm2normn}.
\begin{lem}
\label{norm2normn}
For $f \in \cH$, let $\cA$ be the  event 
\begin{equation}
\label{A}
\cA=\left\{ \forall f \in \cF, \forall v \in \cP, \;
 \|\widehat{f}_v-{f}_v\|_n \leq  2\| \widehat{f}_v-{f}_v\|_2
+\gamma_v\right\}.
\end{equation}
Then, for some positive constant $c_{2}$,
\begin{equation*}
\PXe \left( \cA \right) \geq 1-\sum_v \exp(-n c_{2}\gamma_v^2).
\end{equation*}
\end{lem}
On the set $\cA$, Inequality \eref{borne}  provides that, for all $K>0$
\begin{eqnarray}
\label{intermediaire}
\frac{1}{C} \| \widehat{f}-m\|_n^2 &\leq &
\sum_{v\in  S_f} \left(\mu_v + 2 \gamma_v \| \widehat{f}_v-{f}_v\|_2 + \gamma_v^2
\right) \\
&\leq &\sum_{v\in  S_f} \left(\mu_v + (1+K)\gamma_v^{2}  +\frac{1}{K}\| \widehat{f}_v-{f}_v\|^{2}_2
\right) , \label{intermediaire1}
\\
&\leq &\sum_{v\in  S_f} \left(\mu_v + (1+K)\gamma_v^{2} \right) +\frac{1}{K}\sum_{v \in \mathcal{P}}\| \widehat{f}_v-{f}_v\|^{2}_2
 \nonumber
\\
\label{intermediaire2}
 &\leq & \sum_{v\in  S_f} \left(\mu_v + (1+K)\gamma_v^{2}\right)   
  + \frac{1}{K}\| \sum_{v \in \mathcal{P}}\widehat{f}_v-{f}_v \|^{2}_2
  .
\end{eqnarray}
Inequality~\eref{intermediaire1} uses the inequality  $2ab \leq
\frac{1}{K} a^{2}+ K b^{2}$ for all positive $K$, and
Inequality~\eref{intermediaire2} uses the orthogonality with respect
to $\bL^{2}(\PX)$.

In the following we have to consider several cases, according to the
rankings of
$ \|\sum_{v\in \mathcal{P}}\widehat{f}_v-{f}_v\|_2$, $\| \sum_{v\in \mathcal{P}}\widehat{f}_v-{f}_v\|_n$ and $d(f)$ defined as follows
\begin{equation}d^2(f) = \max\left( \sum_{v\in
    S_f} \gamma_v^2,\sum_{v\in S_f}\mu_v   \right)\label{dn.eq}.
    \end{equation}
More precisely, we consider three cases
\begin{itemize}
\item[\underline{Case 1:}] $ \|  \sum_{v\in \mathcal{P}}\widehat{f}_v-{f}_v\|_2\leq \| \sum_{v\in \mathcal{P}}\widehat{f}_v-{f}_v\|_n.$
\item[\underline{Case 2:}] $\| \sum_{v\in \mathcal{P}}\widehat{f}_v-{f}_v\|_n\leq \| \sum_{v\in \mathcal{P}}\widehat{f}_v-{f}_v\|_2 \leq d(f)$
\item[\underline{Case 3:}] $\| \sum_{v\in \mathcal{P}}\widehat{f}_v-{f}_v\|_n\leq \| \sum_{v\in \mathcal{P}}\widehat{f}_v-{f}_v\|_2 $ and $ d(f) \leq\| \sum_{v\in \mathcal{P}}\widehat{f}_v-{f}_v\|_2$.
\end{itemize}

\medskip
\noindent
\label{3Cases}
\underline{Case 1:}
From \eref{intermediaire2}, for any $f \in  \mathcal{H}$, we get
\begin{eqnarray*}
\frac{1}{C}\| \widehat{f}-m\|_n^2 
& \leq & 
 \sum_{v\in  S_f} \left(\mu_v+ (1+K)
\gamma_v^2 \right)+
  \frac{1}{K}\| \widehat{f}-{f}\|_n^2.
\end{eqnarray*}
Hence, using that for all $K'>0$, 
\begin{equation}
\label{trick2}\|\widehat{f}-f\|_{n}^{2}\leq (1+K') \|\widehat{f} -m\|_{n}^2+
(1+1/K')\|f-m\|_{n}^2,
\end{equation} we obtain for a 
suitable choice of $K'$, say $1+K' < K/C$, that, for some positive constant $C'$,
\begin{eqnarray*}
\| \widehat{f}-m\|_n^2 
  &\leq&  C'\left[
\| f-m\|_n^2+ 
\sum_{v\in  S_f}\mu_v+ 
 \sum_{v\in S_f}\gamma_v ^2\right].
\end{eqnarray*}
This shows the result in Case 1.\\ 

\medskip
\noindent
\underline{Case 2:}
Inequality \eref{intermediaire2} becomes

\begin{equation*}
 \frac{1}{C} \| \widehat{f}-m\|_n^2 \leq 
\sum_{v\in  S_f} \left(\mu_v + (1+K)\gamma_v^{2}\right)   
  + \frac{1}{K} d^2(f),
\end{equation*}
which gives  the expected result since $d^{2}(f)= \max\left\{
  \sum_{v\in  S_f} \mu_v , \sum_{v\in  S_f}\gamma_{v}^{2}\right\}$.

\medskip
\noindent
\underline{Case 3:}
Recall that in this case, $\| \sum_{v\in
  \mathcal{P}}\widehat{f}_v-{f}_v\|_n\leq \| \sum_{v\in
  \mathcal{P}}\widehat{f}_v-{f}_v\|_2 $ and $d(f) \leq\|
\sum_{v\in \mathcal{P}}\widehat{f}_v-{f}_v\|_2$. This case is solved by
applying the following Lemma (proved in
Section~\ref{prooflemcompnormes3},
page~\pageref{prooflemcompnormes3}), which states that with high
probability, $\|\widehat{f} -f \|_{2} \leq \sqrt{2} \|\widehat{f} -f \|_{n}$.
\begin{lem}
\label{lemcompnormes3}
Let $f=\sum_v f_v \in \cF$ with support $S_{f}$, $d(f)$ be
defined by~\eref{dn.eq},  and let
$\mathcal{G}(f)$ be the class of functions written as $g=\sum_{v\in
  \mathcal{P}} g_v$, such that $\| g_v\|_{\mathcal{H}_v}\leq 2$
satisfying for all $f \in \cF$
\begin{eqnarray*}
\mbox{\bf C1}  && \sum_{v\in \mathcal{P}} \mu_v\| g_v \|_{\mathcal{H}_v} +\sum_{v\in \mathcal{P}} \gamma_v\| g_v\|_{n}\leq  
\sum_{v\in S_f} 4\mu_v\| g_v\|_{\mathcal{H}_v}+
\sum_{v\in S_f} 4\gamma_v\| g_v\|_{n} \\
\mbox{\bf C2}  && \sum_{v\in S_f}\gamma_v\| g_v\|_n\leq 2\sum_{v\in
  S_f}\gamma_v\| g_v\|_2+\sum_{v\in S_f} \gamma_v^2\\
\mbox{\bf C3}  && \| g\|_n\leq \| g \|_2.
\end{eqnarray*}
Then the event 
$$\left\{  \| g\|_n^2 \geq \frac{\| g\|^2_2}{2}, \quad \| g\|_2\geq
d(f)
\right\}$$
has probability greater than 
$1-\exp(-n c_3 \sum_{v \in S_{f}} \lambda_{n,v}^2)$ for some constant $c_3$.
\end{lem}
Note that Assumption $n \lambda_{n,v}^{2} \geq -C_{2}\log(\lambda_{n,v})$
implies that  $\lambda_{n,v}= K_{n,v}/\sqrt{n}$ with $K_{n,v}\rightarrow
\infty$. Then,  if $f$ is such that $|S_{f}| \geq 1$, $\exp(-n c_{3}
\sum_{v \in S_{f}} \lambda_{n,v}^2) \leq \exp(-c_{3}\min_{v \in \cP} K_{n,v})$.
If $f$ is such that $|S_{f}|=0$, then Condition {\bf C1} is not
satisfied except if $g_{v}=0$ for all $v  \in \cP$. Because we will
apply Lemma~\ref{lemcompnormes3} to $g_{v}= \widehat{f}_{v} - f_{v}$,
this event has probability 0. Therefore the event 
\begin{equation}
\label{C}
\mathcal{C}=\left\{ \forall f \in \cF, \mbox{ such that } 
 g=\sum_{v\in \mathcal{P}}(\widehat{f}_v-f_v )\in \cG(f),
\mbox{ and }
 \| g\|_n^2 \geq \frac{\| g\|^2_2}{2}, \; \|
 g\|_2^2\geq  d(f)  \right\}
\end{equation}
 has probability greater than $1-\eta/3$ for some $0<\eta<1$.

Conditionning on the events $\mathcal{T}$ and  $\mathcal{A}$
(defined by \eref{evtTau} and \eref{A}),  and according to Lemma \ref{lemme2}, $\sum_{v\in \mathcal{P}}
(\widehat{f}_v-{f}_v)$ belongs to the set $\mathcal{G}(f)$.
According to \eref{intermediaire2}, we conclude in the same way as in
Case 1.

Finally, it remains to quantify $\PXe(\cT \cap \cA  \cap \cC)$. Following
Lemma~\ref{Tau}, and Lemma~\ref{norm2normn}, $\cT$, respectively $\cA$,
has probability greater than $1 - c_1\sum_{v
  \in \mathcal{P}}\exp(-nc_2 \lambda_{n,v}^2)$, respectively
$1-\sum_v \exp(-n \gamma_v^2)$. Each of these probabilities is greater
$1-\eta/3$ thanks to the assumption $n \lambda_{n,v}^{2} \geq -C_{2}
\log \lambda_{n,v}$.

\subsection{\label{PrCorr2.1.st}Proof of Corrolary  \ref{oracle2}}
We start from Theorem \ref{oracle} which states that with high probability,
$$\| \widehat{f}-m\|_n^2 \leq  C\left\{\| f-m\|_n^2+ \sum_{v\in  S_f}\mu_v+\sum_{v\in S_f}\gamma_v ^2 \right\},$$
and  use that for all $\theta>0$,
\begin{eqnarray}
\label{trick}
\| \widehat{f}-m\|_2^2\leq (1+\theta)\| \widehat{f}-f\|_2^2+
(1+\frac{1}{\theta})\| m-f\|_2^2.
\end{eqnarray}
For $d$ defined by \eref{dn.eq}, we consider once again the three
cases defined page~\pageref{3Cases}.

\medskip
\noindent
\underline{Case 1:}
According to  \eref{trick2} and \eref{trick}, we get the result since
$$\| \widehat{f}-m\|_2^2\leq (1+\theta)\parallel \widehat{f}-f\parallel_n^2+
(1+\frac{1}{\theta})\| m-f\|_2^2.
$$

\medskip
\noindent
\underline{Case 2: } We directly obtain that
$$\| \widehat{f}-m\|_2^2\leq (1+\theta)d^2(f)+
(1+\frac{1}{\theta})\| m-f\|_2^2.$$

\medskip
\noindent
\underline{Case 3: }
Recall that in this case, $\| \sum_{v\in \mathcal{P}}\widehat{f}_v-{f}_v\|_n\leq \| \sum_{v\in \mathcal{P}}\widehat{f}_v-{f}_v\|_2 $ and $d(f) \leq\| \sum_{v\in S}\widehat{f}_v-{f}_v\|_2$.
Apply  Lemma \ref{lemcompnormes3} (page~\pageref{lemcompnormes3}) and
conclude that conditionning on the events $\mathcal{T}$ and
$\mathcal{A}$, defined by \eref{evtTau} and \eref{A},
 then $\sum_{v\in \mathcal{P}}\widehat{f}_v-{f}_v$ belongs to $\mathcal{G}(f)$ defined in Lemma \ref{lemcompnormes3}. 
Now, conditionning on the event $\cC$
we get the result since $$\| \sum_{v \in \mathcal{P}} \widehat{f}_v-{f}_v\|_2^2\leq 2\|   \sum_{v \in \mathcal{P}}     \widehat{f}_v-{f}_v\| _n^2.$$
\hfill $\Box$

\subsection{Rate of convergence \label{ratecvg.proof}}
Recall that we consider the case where the
variables $X_{1}, \ldots, X_{d}$ have the same distribution $P_{1}$
on $\cX_{1} \subset\R$, and where the 
unidimensionnal kernels $k_{0 a}$ are all identical.

In this context, our goal is to show that the rate $\nu_{n, v}$ defined at Equation~\eref{nu} is bounded above by a term of
order $n^{-\alpha/(2\alpha + 1)} (\log n)^{\gamma}$, where $\gamma \geq
(|v|-1)\alpha/(2\alpha-1)$.

We start from the fact that, 
$k_{v}(\gx_{v}, \gx'_{v}) = \prod_{a \in v}
k_{0}(x_{a},x'_{a})$, with  a kernel $k_{0}$ admitting an eigen expansion
given   by 
\begin{equation*}
 k_{0}(x,x') = \sum_{\ell \geq 1} \omega_{0, \ell} \zeta_{\ell}(x)\zeta_{\ell}(x') 
\end{equation*}
where the eigenvalues $\omega_{0, \ell}$ are non negative and ranged
in the decreasing order at the rate $\ell^{-2\alpha}$ for some $\alpha > 1/2$, and where the $\zeta_{\ell}$ are the
associated eigen functions, orthonormal with respect to $\bL^{2}(P_1)$.

Therefore the kernel $k_{v}$  admits the following expansion
\begin{equation*}
 k_{v}(\gx_v,\gx^{'}_{v}) = \sum_{\gell=(\ell_1 \ldots \ell_{|v|})} 
\underbrace{\prod_{a=1}^{|v|}\omega_{0, \ell_{a}}}_{\omega_{v, \gell}}
\underbrace{\prod_{a=1}^{|v|}\zeta_{\ell_{a}}(x_{a})}_{\zeta_{v, \gell}(\gx_{v})}
\underbrace{\prod_{a=1}^{|v|}\zeta_{\ell_{a}}(x^{'}_{a})}_{\zeta_{v, \gell}(\gx^{'}_{v})}.
\end{equation*}
According to this expansion the 
$\omega_{v, \gell}$ are of order
$(\prod_{a=1}^{|v|}\ell_{a})^{-2\alpha}$.

In order to control the rate  $\nu_{n,v}$ defined
as $$\nu_{n,v}= \inf\left\{ t \mbox{ such that } Q_{n,v}(t)\leq
  \Delta t^{2}\right\}$$, we have to calculate an upper bound for 
$Q_{n, v}^{2} (t)$
We start with 
the following inequalities that hold up to some constant, for $t^{-1/\alpha}>1$
\begin{eqnarray}
Q_{n, v}^{2} (t) & = & \frac{5}{n} \sum_{\gell} \min(t^{2}, \omega_{v,
  \gell})  \nonumber \\
& \lesssim & \frac{1}{n} t^{2} \sum_{\gell=(\ell_1 \ldots \ell_{|v|})} I (\ell_{1}^{-2\alpha}\times\ldots \times\ell_{|v|}^{-2\alpha} \geq
  t^{2}) + \frac{1}{n} \sum_{\gell=(\ell_1 \ldots \ell_{|v|})} \prod_{a=1}^{|v|}
         \ell_{a}^{-2\alpha} \nonumber \\
& \lesssim & \frac{1}{n} t^{2} \sum_{\gell=(\ell_1 \ldots \ell_{|v|})} I (\ell_{1}\times\ldots \times\ell_{|v|} \leq
  t^{-1/\alpha}) + \frac{1}{n} \left(\sum_{j\geq
         1}\frac{1}{j^{2\alpha}}\right)^{|v|} \label{Q2.eq}
\end{eqnarray}

Now let us mention that $\alpha> 1/2$, $\sum_{j\geq 1}\frac{1}{j^{2\alpha}}$ is a constant
that depends on $\alpha$.  We thus focus on the  first term in the right hand
side of Equation~\eref{Q2.eq}.

Let $u=t^{-1/\alpha} \geq 1$ and let $B_{\vert v \vert}$ be defined as
follows:
$$ B_{\vert v\vert}= \sum_{\gell=(\ell_1 \ldots \ell_{|v|})}  I(\ell_1 \leq u, \ldots,\ell_{|v|} \leq u ) .$$
Let us prove that 
\begin{equation}
B_{\vert v\vert } 
\leq u \left(1+\log(u))^{|v|-1}\right).
\label{Bv.eq}
\end{equation}

Proof of Equation~\eref{Bv.eq}:
First note that
\begin{equation*}
B_{1} = \sum_{\ell\geq 1}  I(\ell \leq u) \leq u I(u\geq 1).
\end{equation*}
In the same way,
\begin{eqnarray*}
B_{2} = \sum_{\ell_{1} \geq 1,  \ell_{2} \geq 1}  
I (\ell_{1}\ell_{2} \leq u) & = & 
\sum_{\ell_{1}\geq 1, \ell_{2}
  \geq 1} 
I(\ell_1 \leq u) I(\ell_2 \leq u/\ell_1) =
 \sum_{\ell_{1}\geq 1}I(\ell_1 \leq u) \sum_{\ell_{2}\geq 1}I(\ell_2 \leq u/\ell_1)\\
&\leq& \sum_{\ell_{1}\geq 1}I(\ell_1 \leq u) \frac{u}{\ell_1}
       I(\frac{u}{\ell_1} \geq 1)
= u \sum_{\ell_{1}\geq 1}\frac{1}{\ell_1} I(\ell_1 \leq u)\\
&\leq& u \left(1 + \log(u)\right).
\end{eqnarray*}
More generally,
\begin{eqnarray*}
B_{\vert v\vert} &= &\sum_{\ell_{1}\geq 1,  \ldots \ell_{|v|} \geq 1}  
 I (\ell_{1} \ldots \ell_{|v|} \leq u) \\
& = & \sum_{\ell_{1}\geq 1, \ldots, \ell_{|v|-1}\geq 1} 
I(\ell_1 \ldots \ell_{|v|-1}\leq u) 
\sum_{\ell_{|v|} \geq 1}I(\ell_{|v|} \leq u/\ell_1... \ell_{|v|-1})\\
&\leq& u \sum_{\ell_{1}\geq 1, \ldots, \ell_{|v|-1}\geq 1} 
\frac{1}{\ell_1 \ldots \ell_{|v|-1}}I(\ell_1 \ldots \ell_{|v|-1}\leq  u)
\end{eqnarray*}

Let $A_{|v|}$ be defined as follows:
$$  A_{\vert v\vert}=\sum_{\ell_{1}\geq 1,  \ldots \ell_{|v|} \geq 1}  
\frac{1}{\ell_1 \ldots \ell_{|v|}} I (\ell_{1} \ldots \ell_{|v|} \leq u), $$
then we get $B_{\vert v\vert} \leq  u A_{\vert v\vert - 1}$. If we show that 
\begin{equation}
 A_{\vert v\vert} 
\leq \left(1 + \log(u)\right)^{|v|}, \label{Av.eq}
\end{equation}
then Inequality~\eref{Bv.eq} is proved.

Proof of Equation~\eref{Av.eq} :
\begin{equation*}
A_1 =  \sum_{\ell\geq 1} \frac{1}{\ell}  I(\ell \leq u) = 1 + \sum_{\ell
   \geq 2} \frac{1}{\ell} I(\ell \leq u) \leq 1 + \int_{1}^{u} \frac{1}{v}
 dv = 1 + \log(u)
\end{equation*}
\begin{eqnarray*}
A_2 =  \sum_{\ell_1 \geq 1, \ell_2 \geq 1} 
\frac{1}{\ell_1 \ell_2}  I(\ell_1 \ell_2 \leq u) & = &
\sum_{\ell_1 \geq 1} \frac{1}{\ell_1} I(\ell_1 \leq u)
\sum_{\ell_2 \geq 1} \frac{1}{\ell_2} I(\ell_2 \leq u/\ell_1) \\
& = &\sum_{\ell_1 \geq 1} \frac{1}{\ell_1} I(\ell_1 \leq u)
\left(1 + \log(u/\ell_1)\right)\\
& = &
\left(1 + \log(u)\right)^{2} - \sum_{\ell_1 \geq 1}
\frac{\log(\ell_1)}{\ell_1}  I(\ell_1 \leq u)\\
& \leq  &\left(1 + \log(u)\right)^{2} .
\end{eqnarray*}
In the same way we have 
\begin{eqnarray*}
A_{\vert v\vert} &= &\sum_{\ell_{1}\geq 1,  \ldots \ell_{|v|} \geq 1}  
\frac{1}{\ell_1 \ldots \ell_{|v|}} I (\ell_{1} \ldots \ell_{|v|} \leq u) \\
& = & \sum_{\ell_{1}\geq 1, \ldots, \ell_{|v|-1}\geq 1} 
\frac{1}{\ell_1 \ldots \ell_{|v|-1}}I(\ell_1 \ldots \ell_{|v|-1}\leq u) 
\sum_{\ell_{|v|} \geq 1}\frac{1}{\ell_{|v|}}I(\ell_{|v|} \leq u/\ell_1... \ell_{|v|-1})\\
&\leq& \sum_{\ell_{1}\geq 1, \ldots, \ell_{|v|-1}\geq 1} 
\frac{1}{\ell_1 \ldots \ell_{|v|-1}}I(\ell_1 \ldots \ell_{|v|-1}\leq  u)
\left(1+\log(u/\ell_1 \ldots \ell_{|v|-1})\right)\\
&\leq& \left(1 + \log(u)\right)^{|v|}.
\end{eqnarray*}
And Bound \eref{Av.eq} is proved.

\paragraph{{\bf Rate of convergence}}
Let us come back to the control of the rate $\nu_{n,v}= \inf\left\{ t \mbox{ such that } Q_{n,v}(t)\leq  \Delta t^{2}\right\}$.
Thanks to~\eref{Bv.eq}  we obtain that, 
up to some constant that depends on $|v|$ and $\alpha$,
\begin{eqnarray*}
 Q_{n,v}^{2} (t) &  \lesssim & \frac{1}{n} t^{2} \sum_{\gell} I (\ell_{1}\times\ldots \times\ell_{|v|} \leq
  t^{-1/\alpha}) + \frac{1}{n} \left(\sum_{j\geq
         1}\frac{1}{j^{2\alpha}}\right)^{|v|}\\
 &  \lesssim & \frac{1}{n} t^{2-1/\alpha} \left(1 - \frac{1}{\alpha}
          \log(t)\right)^{|v|-1} + \frac{1}{n}
\end{eqnarray*}
It remains now to find $t$ such that, up to constant
 \begin{equation*}
 \frac{1}{\sqrt{n}} t^{1-1/2\alpha} 
\left(1 - \frac{1}{\alpha} \log(t)\right)^{(\vert v\vert -1)/2}
\leq
 t^{2}
\end{equation*}

If $t = n^{-\beta} (\log(n))^{\gamma}$ with $\beta =
\alpha /(1+2\alpha)$, $\gamma >0$, $\alpha>1/2$, then
\begin{eqnarray*}
 1 - \frac{1}{\alpha} \log(t) &= &
1 - \frac{1}{\alpha} \log\left(n^{-\alpha /(1+2\alpha)}
                                   (\log(n))^{\gamma}\right) \\
&= &
1 - \frac{1}{\alpha}\left(
-\frac{\alpha}{1+2\alpha}\log(n) +\gamma \log\log(n)
\right)\\
&= &
1 + \frac{1}{1+2\alpha}\log(n) -\frac{\gamma}{\alpha} \log\log(n)
\\&\leq  &\log(n) \mbox{ as soon as } \log(n) > 1 + \frac{1}{2\alpha}.
\end{eqnarray*}

Therefore $\nu_{n,v}$  will be smaller than the infimum of $t$ such
that 
\begin{equation*}
 \frac{1}{\sqrt{n}} t^{1-1/2\alpha} \left(\log(n)\right)^{(\vert
   v\vert -1)/2} \leq t^{2},
\end{equation*}

which is satisfied if  $\gamma \geq (|v|-1) \alpha/(2\alpha -1)$.

\subsection{Intermediate Lemmas}
\label{intermediate}
For $v\in \mathcal{P}$, let $\cH_{v}$ be the RKHS associated to the
self reproducing kernel 
$k_v$. Let $Q_{n,v}$ and $\nu_{n, v}$ and  be defined by
Equations~\eref{Qn} and~\eref{nu}. For any function $g_v \in \cH_{v}$, let $V_{n, \varepsilon}$ be defined at Equation~\eref{procemp} and consider the following processes
\begin{eqnarray}
 W_{n,2,v}( t) & = & \sup \left\{ \left|V_{n, \varepsilon}(g_{v})\right|, \;
\|g_v\|_{\mathcal{H}_v}\leq 2, \|g_v\|_{2}\leq t\right\} \label{defWn} \\
W_{n,n,v}( t) & =& \sup \left\{ \left|V_{n, \varepsilon}(g_{v})\right|, 
\|g_v\|_{\mathcal{H}_v}\leq 2, \; \|g_v\|_{n}\leq t\right\}.\label{defWntilde}
\end{eqnarray}

\begin{lem}
\label{lemcomplex}
If $\EXe$ denotes the expectation with
respect to the distribution of $(\gX, \varepsilon)$,  we have
for all $t > 0 $,
\begin{equation*}
\EXe  W_{n,2,v}(t)
 \leq  Q_{n,v}(t).
\end{equation*}
\end{lem}
Its proof is given in Section~\ref{Prooflemcomplex}
page~\pageref{Prooflemcomplex}.

\begin{lem}
\label{lemcompnormes1}
Let  $b > 0$ and let $\cG(t)$ be the following class of functions:
\begin{equation}
\cG(t)=\left\{ g_v \in \cH_{v}, \|g_v\|_{\cH_v}\leq 2, \|g_v\|_{2}\leq t,
  \|g_v\|_{\infty}\leq b \right\}.
\label{calG.eq}
\end{equation}
Let $\Omega_{v,t}$ be the event defined as
\begin{eqnarray}
\label{Omega}
\Omega_{v,t}= \left\{
\sup \left\{
\vert \| g_v\|_2- \| g_v\|_n \vert, \; g_v \in \cG(t) \right\} \leq \frac{bt}{2}\right\}.
\end{eqnarray}
Then for any $t\geq \nu_{n,v}$, the event $\Omega_{v,t}$ has
probability greater than $ 1-\exp(-c_2 n t^2)$, 
for some positive constant $c_{2}$.
\end{lem}
Its proof is given in Section~\ref{Prooflemcompnormes1},
page~\pageref{Prooflemcompnormes1}.

\begin{lem}
\label{lemcompnormes2}
For any function $g_v \in \mathcal{H}_v$
satisfying $\| g_v\|_{\mathcal{H}_v}\leq 2$,   $\| g_v\|_\infty \leq b$ and 
 $\| g_v\|_2\geq t$, for all $t\geq
   \nu_{n,v}$ and $b \geq 1$, the event
\begin{equation*}
\left(1 - \frac{b}{2}\right) \| g_v\|_2 \leq \| g_v\|_n \leq 
\left(1 + \frac{b}{2}\right) \| g_v\|_2
\end{equation*}
has probabilty greater than $1-\exp(-c_2 n t^2)$ for some positive
constant $c_{2}$.
\end{lem}
Its proof is given in Section~\ref{Prooflemcompnormes2},
page~\pageref{Prooflemcompnormes2}.

\begin{lem}
\label{concentration1}
If $\Ee$ denotes the expectation with
respect to the distribution of $\varepsilon$,  we have
\begin{eqnarray}
\label{concentnn}\PXe \left\{ 
\left\vert W_{n,n,v}(t) - \Ee \big( W_{n,n,v}(t)\big) \right\vert\geq
\delta t   \right\} \leq 4\exp\left(- \frac{n\delta^2}{2}\right).
\end{eqnarray}
\end{lem}
Its proof is given in Section~\ref{Proofconcentration1},
page~\pageref{Proofconcentration1}.

\begin{lem}
\label{concentration2}
Conditionnaly on the space $\Omega_{v,t}$ defined by \eref{Omega}, we
have the two following inequalities:
\begin{eqnarray}
\label{concentn2}
\PXe \left\{ \left\vert W_{n,2,v}(t) - \Ee\big( W_{n,2,v}(t)\big) 
\right\vert \geq \delta t   \right\} & \leq & 4\exp\left(-
  \frac{n\delta^2}{8}\right) \\
\label{concentration}
\PX \left\{  \Ee W_{n,2,v}( t) - \EXe
\big( W_{n,2,v}(t)\big)  \geq  x \right\} &\leq &
\exp\left( -\frac{n x^2}{Q_{n,v}(t)} \right).
\end{eqnarray}
\end{lem}
Its proof is given in Section~\ref{Proofconcentration2},
page~\pageref{Proofconcentration2}.

 \begin{lem}
  \label{Case1} Let $\lambda_{n,v}$ be defined at
  Equation~\eref{lambda.eq}, $\Delta$ at Equation~\eref{nu} and
  $\kappa= 10+4\Delta$. Conditionnaly
on the space  $\Omega_{v,\lambda_{n,v}}$ defined at Equation~\eref{Omega}, for some positive constants $c_1, c_{2}$, 
 with probability greater than $1-c_{1}\exp(-c_{2} n\lambda_{n,v}^2)$,
we have 
\begin{equation}
 \label{etap} W_{n,n,v}( \lambda_{n,v}) \leq \kappa \lambda_{n,v}^{2}\;
 \mbox{ and } \; \Ee W_{n,n,v}( \lambda_{n,v}) \leq \kappa \lambda_{n,v}^{2}.
\end{equation}
\end{lem}
Its proof is given in Section~\ref{ProofCase1},
page~\pageref{ProofCase1}.

\subsection{Proofs of  Lemma \ref{Tau} to
~\ref{lemcompnormes3}: }
\label{proofsintermediate}

\subsubsection{\label{ProofTau}Proof of Lemma \ref{Tau} (page~\pageref{Tau})}

For $f \in \cF$ and $v \in \cP$, let $g_{v} =
\widehat{f}_v-f_v$. Note that $\|g_{v}\|_{\cH_{v}} \leq 2$. Let us show that 
\begin{equation}
\label{but}
\vert V_{n,\varepsilon}(g_v)\vert 
\leq 
\kappa\left[\lambda_{n,v}^2\| g_v\|_{\mathcal{H}_v}+\lambda_{n,v}\| g_v\|_{n}\right].
\end{equation}

We start by writing that
\begin{equation}
\left| V_{n,\varepsilon}(g_v) \right|=
\| g_v\|_{\mathcal{H}_v}\left| V_{n,\varepsilon}\left(\frac{g_v}{\| g_v\|_{\mathcal{H}_v}}\right) \right|
\leq \| g_v\|_{\mathcal{H}_v} W_{n,n,v}\Big(\frac{\| g_v\| _{n}}{ \|
  g_v\|_{\mathcal{H}_v}}\Big).
\label{baseW}
\end{equation}
Consider the two following cases:
\begin{itemize}
\item[\underline{Case A:}] $\| g_v\|_n \leq \lambda_{n,v} \|
  g_v\|_{\mathcal{H}_v}$
\item[\underline{Case B:}] $\| g_v\|_n > \lambda_{n,v} \| g_v\|_{\mathcal{H}_v}$
\end{itemize}

\medskip
\noindent
\underline{Case A}: \label{CaseA}
Since $\| g_v\| _{n}\leq \lambda_{n,v}\| g_v\|_{\mathcal{H}_v} $, we have
$$W_{n,n,v}\Big(\frac{\| g_v\| _{n}}{ \|
  g_v\|_{\mathcal{H}_v}}\Big)\leq W_{n,n,v}( \lambda_{n,v}).$$ We
then apply Lemma \ref{Case1}, page~\pageref{Case1}, and conclude that \eref{but} holds in Case A for each  $v\in \mathcal{P}$ since,
 with high probability
\begin{eqnarray}
\label{cas1}
\vert V_{n,\varepsilon}(g_v)\vert 
&\leq& 
\kappa\lambda_{n,v}^2\| g_v\|_{\mathcal{H}_v}\leq 
\kappa\lambda_{n,v}^2\| g_v\|_{\mathcal{H}_v}+
\kappa\lambda_{n,v}\| g_v\|_{n}.
  \end{eqnarray}

\medskip
\noindent\underline{Case B}: Consider now the case $\| g_v\| _{n}>\lambda_{n,v}
\| g_v\| _{\mathcal{H}_v}$ and let us show that for any $v\in \mathcal{P}$,
\begin{eqnarray*}
W_{n,n,v}\left(\frac{\| g_v\|_n}{\| g_v\|_{\mathcal{H}_v}}\right)
 \leq 
\kappa\lambda_{n,v}\| g_v\|_{n}
 .
\end{eqnarray*}
Let $r_v$ be a deterministic number such that $r_v>\lambda_{n,v}$. Our first step relies on the study of the process
 $W_{n,n,v}(r_v),$ for  $r_v>\lambda_{n,v}$.
In that case we state two results:
\begin{itemize}
\item[{\bf R1}]  For any deterministic $r_v\geq \lambda_{n,v}$, with probability greater than $1-c_1 \exp(-c_2 n \lambda_{n,v}^2)$,
\begin{equation}\label{int43}
W_{n,n,v}(r_v) \leq \kappa r_v \lambda_{n,v}.
\end{equation}

\item[{\bf R2}] Inequality~\eref{int43} continues to hold
  for random $r_v$ of the form $$r_v=\frac{\| g_v\|_n}{\|
    g_v\|_{\mathcal{H}_v}}.$$
\end{itemize}
Combining these two points implies that, with probability greater than $1-c_1 \exp(-c_2 n \lambda_{n,v}^2)$,
\begin{equation*}
 \| g_v\|_{\mathcal{H}_v}W_{n,n,v}\Big(\frac{\| g_v\|_n}{\|
  g_v\|_{\mathcal{H}_v}}\Big) \leq
\kappa \| g_v\|_n \lambda_{n,v}.
\end{equation*}

Consequently, in Case B, according to \eref{baseW}, for each $v$,
Inequality~\eref{but} holds because
\begin{equation*}
\vert V_{n,\varepsilon}(g_v)\vert 
\leq 
\kappa \| g_v\|_n \lambda_{n,v}
\leq 
\kappa \lambda_{n,v}^2\| g_v\|_{\mathcal{H}_v}+ 
\kappa \lambda_{n,v}\| g_v\|_{n}.
\end{equation*}
 This ends up the proof of Lemma \ref{Tau}.

\medskip
\paragraph{\underline{Proof of {\bf R1}}}
Taking $t=r_v$ and $\delta=\lambda_{n,v}$ in  \eref{concentnn}, with probability greater than 
 $1-4\exp(-n\lambda_{n,v}^2)$, we
have $$W_{n,n,v}(r_v) \leq    \Ee [W_{n,n,v}(r_v)]+ r_v \lambda_{n,v}.$$
  
 Next we prove that for some positive  $r_v$, with probability greater
 than $1-\exp(-nc\lambda_{n,v}^2)$, we have
 \begin{eqnarray}
\label{intcas2}\Ee W_{n,n,v}( r_v)  \leq
\kappa  r_v \lambda_{n,v}.
\end{eqnarray}
Let $\widehat{\nu}_{n,v}$ defined 
as the smallest solution of $\Ee [ W_{n,n,v}(t)]\leq \kappa t^2$.
 For $W_{n,n,v}$, defined by \eref{defWn}, we write
\begin{eqnarray*}
\Ee W_{n,n,v}(r_v)
&=& \frac{r_v}{\widehat{\nu}_{n,v}}
\Ee  \sup \left\{
\vert V_{n,\varepsilon}(g_v)\vert , \;
\| g_v\|_{\mathcal{H}_v}\leq \widehat{\nu}_{n,v}/r_v , \;
\| g_v\|_{n}\leq \widehat{\nu}_{n,v}
\right\}
\\
&\leq& \frac{r_v}{\widehat{\nu}_{n,v}}  \Ee
 W_{n,n,v}(\widehat{\nu}_{n,v}) \leq \frac{r_v}{\widehat{\nu}_{n,v}}
\kappa \widehat{\nu}_{n,v} ^2
=
\kappa  r_v \widehat{\nu}_{n,v}.
\end{eqnarray*}
Besides, Lemma~\ref{Case1} stated that on the event $\Omega_{v,
  \lambda_{n,v}}$, $\Ee W_{n,n,v}( \lambda_{n,v}) \leq \kappa
\lambda_{n,v}^{2}$. It follows from the definition of
$\widehat{\nu}_{n,v}$, and Lemma~\ref{lemcompnormes1}, that $\widehat{\nu}_{n,v} \leq \lambda_{n,v}$
for all $v \in
\mathcal{P}$ with probability greater than $1 - \exp(-n c_{2} \sum_{v
  \in \cP} \lambda_{n,v}^{2})$. 
 Consequently, for any deterministic $r_v$ such that $r_v\geq
 \lambda_{n,v}$, \eref{int43} is satisfied with high probability.

\medskip
\paragraph{\underline{Proof of {\bf R2}}}
Let us  prove  {\bf R2} by using a peeling-type argument.  Our aim is to prove that \eref{int43} holds for any $r_v$ of the form
$$r_v=\frac{\| g_v\|_n}{\| g_v\|_{\mathcal{H}_v}}.$$ 
Since $\|g_v\|_\infty/\| g_v\|_{\mathcal{H}_v} \leq 1$, we have
$\|g_v\|_n/\| g_v\|_{\mathcal{H}_v} \leq 1$. 
We thus restrict ourselves  to $r_v$ satisfying $r_v= \| g_v\|_n/\|
g_v\|_{\mathcal{H}_v}$ with $\| g_v\|_n/\| g_v\|_{\mathcal{H}_v} \in (\lambda_{n,v},1]$.

We start by splitting the interval $(\lambda_{n,v},1]$ into $M$
disjoint intervals such that 
$(\lambda_{n,v},1]=\cup_{k=1}^M(2^{k-1}\lambda_{n,v},2^k\lambda_{n,v}]$, for some
$M$ that will be chosen later.
Consider the event $\cD^{c}$ defined as follows:
$$\mathcal{D}^c=\{ \exists v \in \mathcal{P} \mbox{ and } \exists  \overline{g}_v, \mbox{ such that }   
\vert V_{n,\varepsilon}(\overline{g}_v)\vert 
\geq  
\kappa \lambda_{n,v}\| \overline{g}_v \|_{n}, \mbox{ with } \frac{\| \overline{g}_v\|_n}{\| \overline{g}_v\|_{\mathcal{H}}} \in (\lambda_{n,v},1]\}.$$

We prove that, for some positive constants $c_1, c_2$,
$$P [\mathcal{D}^c] \leq c_1\exp(-c_2 n\lambda_{n,v}^2).$$

For $\overline{g}_v \in \cD^{c}$, let $\overline{k}$ be the integer in $\{1,\cdots,M\}$, such that
$$2^{\overline{k}-1}\lambda_{n,v} \leq   \frac{\| \overline{g}_v\|_n}{\| \overline{g}_v\|_{\mathcal{H_{v}}}}  \leq 2^{\overline{k}}\lambda_{n,v}.$$
This $\overline{k}$ satisfies
\begin{eqnarray*}
\| \overline{g}_v\|_{\mathcal{H}_v}W_{n,n,v}\Big(  2^{\overline{k}}\lambda_{n,v}\Big)&\geq&
\| \overline{g}_v\|_{\mathcal{H}_v}W_{n,n,v}\Big(  \frac{\| g_v\|_n}{\| g_v\|_{\mathcal{H}_v}}\Big)
\geq \vert V_{n,\varepsilon}(\overline{g}_v)\vert 
\geq
\kappa \lambda_{n,v}\| \overline{g}_v\|_{n}.
\end{eqnarray*}
Therefore, we get
$$W_{n,n,v}(2^{\overline{k}}\lambda_{n,v}) \geq
\kappa \lambda_{n,v}\frac{\|
  \overline{g}_v\|_n}{\| \overline{g}_v\|_{\mathcal{H}_v}}\geq
\kappa \lambda_{n,v}^2 2^{\overline{k}-1}
\geq 
\kappa \frac{\lambda_{n,v}}{2} 2^{\overline{k}}\lambda_{n,v}.$$
By taking $r_v=2^{\overline{k}}\lambda_{n,v}$ in \eref{int43}, we have
$$\mathcal{P} \left[ W_{n,n,v}(2^{\overline{k}}\lambda_{n,v}) 
\geq 
\kappa \frac{\lambda_{n,v}}{2}
2^{\overline{k}}\lambda_{n,v}\right]\leq c_1 \exp(-c_2n \lambda_{n,v}^2).$$

Now let us write $\mathcal{D}^c$ as
follows:
$$ \mathcal{D}^c =\bigcup_{k=1,\cdots,M}\left\{\; 
\exists v, \exists \; \overline{g}_v \mbox{ such that}
\vert V_{n,\varepsilon}(\overline{g}_v)\vert 
\geq 
\kappa \lambda_{n,v}\frac{\| \overline{g}_v\|_n}{\| \overline{g}_v\|_{\mathcal{H}_v}}
\mbox{ with } \frac{\| \overline{g}_v\|_n}
{\| \overline{g}_v \|_{\mathcal{H}}} 
\in (2^{k-1}\lambda_{n,v},2^k\lambda_{n,v})\right\}
.$$
The set
$\mathcal{D}^c$ has probability smaller than  $c_1 M \exp(-c_2n \lambda_{n,v}^2)$.
If we choose $M $ such that $\log M \leq (c_{2}/2) n
\lambda_{n,v}^{2}$, then the probability of the set $\cT$ is  greater than 
\begin{equation*}
1-\sum_{v \in\cP}  c_{1} \exp- \frac{c_{2}}{2} n \lambda_{n,v}^{2}.
\end{equation*}

It follows that \textbf{R2} is proved which ends up the proof of Lemma \ref{Tau}.
\hfill $\Box$

\subsubsection{\label{Prooflemme2}Proof of Lemma  \ref{lemme2} (page~\pageref{lemme2}).}
Starting from 
\eref{base} with $B$ defined by Equation~\eref{B.eq}, we write
\begin{eqnarray*}
 \frac12\| \widehat{f}-m\|_n^2 & \leq &
2\vert V_{n,\varepsilon}(\widehat{f}-fv)\vert + \\
&& \sum_{v\in S_f}[\mu_v\|\widehat{f}_v-f_v\|_{\mathcal{H}_v}+\gamma_v \|\widehat{f}_v-f_v\|_{n}]- \sum_{v\in S^c}
[\mu_v\| \widehat{ f}_v\|_{\mathcal{H}_v} +\gamma_v \| \widehat{f}_v\|_{n}].
\end{eqnarray*}
On the event $\mathcal{T}$ defined in \eref{evtTau} we have
\begin{eqnarray*}
\frac12 \| \widehat{f}-m\|_n^2 &\leq&  
2 
\kappa \sum_{v\in \mathcal{P}} \lambda_{n,v}^2 \|
\widehat{f}_v-f_v\|_{\mathcal{H}_v}+2
\kappa\sum_{j\in \mathcal{P}} \lambda_{n,v}\| \widehat{f}_v-f_v\|_{n}+\\
&&\sum_{v\in S_f}[\mu_v  \| \widehat{f}_v-f_v        \|_{\mathcal{H}_v}+\gamma_v \| \widehat{f}_v-f_v       \|_{n}]
- \sum_{v\in S^c}
[\mu_v\| \widehat{ f}_v\|_{\mathcal{H}_v}+\gamma_v\| \widehat{ f}_v\|_{n}].
\end{eqnarray*}
Rearranging the terms we obtain that
\begin{eqnarray*}
\frac12\| \widehat{f}-m\|_n^2 &\leq& 
\sum_{v\in S_f}(2 
\kappa\lambda_{n,v}^2+\mu_v)\|
\widehat{f}_v-f_v\|_{\mathcal{H}_v}+\sum_{v\in S_f}(2
\kappa\lambda_{n,v}+\gamma_v)\| \widehat{f}_v-f_v\|_{n}+\\
&&\sum_{v\in S_{f}^c}
(2 
\kappa\lambda_{n,v}^2 -\mu_v)\| \widehat{ f}_v\|_{\mathcal{H}_v}+\sum_{v\in S_{f}^c}
(2 
\kappa\lambda_{n,v} -\gamma_v)\| \widehat{ f}_v\|_{n}.
\end{eqnarray*}
Now, thanks to Assumption~\eref{condmu} with $C_{1} \geq \kappa$
  we have  $ 
\kappa\lambda_{n,v}^2
  \leq \mu_v$  and $2 
\kappa\lambda_{n,v}
  \leq\gamma_v$ and Lemma \ref{lemme2} is shown  since
\begin{eqnarray*}
0\leq \frac12\| \widehat{f}-m\|_n^2 &\leq & 
3\sum_{v\in S_f}
\mu_v\| \widehat{f}_v-f_v\|_{\mathcal{H}_v}
+3\sum_{v\in S_f}
\| \widehat{f}_v-f_v\|_{n}-
\\ &&\sum_{v\in S_{f}^c} \mu_v\| \widehat{ f}_v\|_{\mathcal{H}_v} -\sum_{v\in S_{f}^c} \gamma_v\| \widehat{ f}_v\|_{n}.
\end{eqnarray*}

\hfill $\Box$

\subsubsection{\label{Proofnorm2normn}Proof of  Lemma \ref{norm2normn}
  (page~\pageref{norm2normn}):}

Let us consider the following two cases:
\begin{itemize}
\item $\|\widehat{f}_v-{f}_v\|_2 \leq \gamma_v$. We apply
  Lemma~\ref{lemcompnormes1} (page~\pageref{lemcompnormes1}) to the
  function  $g_{v}=\widehat{f}_v-{f}_v$. It satisfies $g_{v}\in
  \cG(\gamma_v)$ with $b=2$ (recall that $\|\cdot\|_{\infty} \leq
  \|\cdot\|_{\cH_{v}}$). Moreover, $\gamma_{v} \geq C_{1}
  \lambda_{n,v}\geq C_{1} \nu_{nv} \geq \nu_{n,v}$ as soon as $C1 \geq
  1$. 

It follows that, for some positive $c_{2}$, with probability greater than $1-\exp(-n c_{2}\gamma_{v}^{2})$
\begin{equation*}
\| \widehat{f}_v-{f}_v\|_n\leq
  \| \widehat{f}_v-{f}_v\|_2+\gamma_v
\end{equation*}

\item $\|\widehat{f}_v-{f}_v\|_2 \geq \gamma_v$. We apply
  Lemma~\ref{lemcompnormes2} (page~\pageref{lemcompnormes2}) to the
  function  $g_{v}=\widehat{f}_v-{f}_v$ with $b=2$. 
It follows that, for
some positive $c_{2}$, with probability greater than $1-\exp(-n c_{2}\gamma_{v}^{2})$,
\begin{equation*}
\| \widehat{f}_v-{f}_v\|_n\leq
  2 \| \widehat{f}_v-{f}_v\|_2.
\end{equation*}
\end{itemize}
\hfill $\Box$

\subsubsection{\label{prooflemcompnormes3}Proof of Lemma
  \ref{lemcompnormes3} (page~\pageref{lemcompnormes3}): }
Let $d(f)$ be defined by \eref{dn.eq}, and let $\mathcal{G}(f)$
and $\mathcal{G}'(f)$ be the following sets
\begin{eqnarray*}
 \mathcal{G}(f) &=&
\left\{ g=\sum_{v\in \mathcal{P}}g_v, \mbox{ satisfying }
 \|g_{v}\|_{\cH_v} \leq 2, \mbox{ and Conditions \textbf{C1, C2, C3}
 }\right\},\\
\mathcal{G}'(f) & = & \left\{ g\in \mathcal{G}(f),
\mbox{ such that } \| g\|_2=d(f)\right\}.
\end{eqnarray*}
Let us consider the two events $\cB$ and $\cB'$ defined as follows:
$$\mathcal{B}'=\{\forall h \in \mathcal{G}', \; \| h\|^{2}_n \geq
d(f)^{2}/2  \}\mbox{ and }
\mathcal{B}=\{\forall h \in \mathcal{G}, \; \| h\|_n^2\geq \| h\|_2^2/2 , \mbox{ and } \| h\|_2\geq d(f)\}.$$
Let us first remark that $\mathcal{B}'$ is included into
$\mathcal{B}$: if $h \in \cB'$ , then $h \in \cG$, $\|h\|_{2} =
d(f)$ and $\| h\|_n^2\geq d^(f)2/2$. It follows that  $\| h\|_n^2\geq
\| h\|_2^2/2$ and $\|h\|_{2} \geq d(f)$.
Therefore Lemma~\ref{lemcompnormes3} is proved if $\mathcal{B'}$ holds
with high probability.
Consider
 $$Z_n( \mathcal{G}')  =\sup_{g\in \mathcal{G}'}\{d(f)^2-\| g\|_n^2\}.$$ 
We show that the event $Z_n( \mathcal{G}')\leq d(f)^2/2$ has
probability greater than $1-c_1  \exp(-nc_2 d(f)^2)$. 

Let us briefly recall the notion of covering numbers for a totally bounded metric space $(\mathbb{G},\rho)$, consisting of a set $\mathbb{G}$ and a metric $\rho$ defined from $\mathbb{G}\times\mathbb{G} $ into $\mathbb{R}_+$. A $\delta$-covering set of $\mathbb{G}$ is a collection of functions $f^1,\cdots,f^N$ such that for all $f\in \mathbb{G}$ there exists $k\in \{1,2,\cdots,N\}$ such that $\rho(f,f^k)\leq \delta$.  

The $\delta$-covering number $N(\delta,\bG,\rho)$ is the
cardinality of the smallest $\delta$-covering set. A propercovering
restricts  the
covering to use only elements in the set $\bG$. The proper covering
number denoted $N_{\rm{pr}}(\delta,\bG,\rho)$
satisfies
\begin{eqnarray}
\label{encadreN}N(\delta,\bG,\rho) \leq N_{\rm{pr}}(\delta,\bG,\rho) \leq N(\delta/2,\bG,\rho). 
\end{eqnarray}
 
Let us now consider a $d(f)/8$-covering  of $(\mathcal{G}', \|\cdot \|_{n})$, so that, for all 
$g$ in $\mathcal{G}'$ there exists $g^k$ such that $\| g-g^k\|_n \leq
d(f)/8.$  The associated proper covering number  is
\begin{equation}
 N_{\rm{pr}}=N_{\rm{pr}}(d(f)/8,\mathcal{G}', \|\cdot \|_{n} ).
\label{CovN.eq}
\end{equation}
Now, for all $g \in \mathcal{G}'$, $
T_1= \| g^k\|_n^2-\| g\|_n^2,$
and $T_2=d^2(f)-\| g^k\|_n^2$,
we write
\begin{equation*}
 d(f)^2-\| g\|_n^2=T_1+T_2.
\end{equation*}

The proof  is splitted into four steps:
\begin{enumerate}
\item[\underline{Step 1}] The first step consists in showing that 
\begin{eqnarray}
\label{etape1}
T_1=\| g^k\|_n^2-\|  g\|_n^2\leq \frac{d(f)^2}{4}.
\end{eqnarray}
\item[\underline{Step 2}] The second step consists in proving that, for $N_{\rm{pr}}$ given at Equation~\eref{CovN.eq} and for some constant $C$,
\begin{eqnarray*}
\PX \left[ \max_{k\in\{1,\cdots,N_{\rm{pr}}\}}[d(f)^2-\| g^k\|_n^2]\geq
  d^2/4   \right]
\leq \exp\Big(\log{N_{\rm{pr}}}- C n d(f)^2 \Big).
\end{eqnarray*}
\item[\underline{Step 3}]
The third  step concerns the control of  $N_{\rm{pr}}$: we show the
following result
\begin{eqnarray*}
\log N_{\rm{pr}}
 \leq n \left( \frac{64}{d(f)}\mathbb{E}_\varepsilon \sup_{g\in \mathcal{G}'} \vert V_{n,\varepsilon}(g)\vert \right)^{2}  .
\end{eqnarray*}
\item[\underline{Step 4}]
The last step consists  in bounding from above the Gaussian
complexity:  
$$\mathbb{E}_\varepsilon  \sup_{g\in \mathcal{G}'}
   \vert V_{n,\varepsilon}(g)\vert \leq
\frac{20 \kappa}{C_{1}} d(f)^2.$$ 
\end{enumerate}

Let us conclude the proof of the lemma before proving  these four steps.

Putting together Steps 3 and 4, for $c_{3} < C$ and $C_{1}$ large
enough, 
then Step 2 states
that
\begin{equation*}
 \PX\left( T_{2} \geq \frac{d(f)^{2}}{4}\right) \leq \PX \left[ \max_{k\in\{1,\cdots,N_{\rm{pr}}\}}[d(f)^2-\| g^k\|_n^2]\geq
  d(f)^2/4   \right]
\leq \exp\Big(- c_{3} n d(f)^{2}\Big).
\end{equation*}
Now, we have
\begin{equation*}
 \PX \left[Z_n( \mathcal{G}')\leq d(f)^2/2\right] =
\PX \left[\max_{g^1,\cdots,g^N} \{  d(f)^2-\| g^k\|_n^2  \}\geq \frac{d(f)^2}{4}\right]\leq \exp\Big(-c_{3} nd(f)^2\Big).
\end{equation*} 
We conclude the proof of the lemma, by noting that $d(f)^{2}  \geq C_{1}^{2}
\sum_{v} \lambda_{n,v}^{2}$
(see~\eref{lambda.eq}, \eref{dn.eq} and  \eref{condmu}).

\medskip
\noindent
\underline{{Proof of Step 1: }}\label{ProofS1}
We start by writing that
\begin{eqnarray*}
\| g^k\|_n^2-\| g\|_n^2&=&\frac{1}{n}
\sum_{i=1}^n \left[ (g^k(\gX_{i}))^2-(g(\gX_{i}))^{2} \right]\\
&\leq& \| g^k- g\|_n \sqrt{\frac{1}{n}\sum_{i=1}^n [g^k(\gX_{i})+g(\gX_{i})]^2}.
\end{eqnarray*}
Using that $(a+b)^2\leq 2a^2+2b^2$, $g \in \cG'$, and $g$ satisfies
Condition {\bf C3}, we get
\begin{eqnarray*}
\frac{1}{n}\sum_{i=1}^n [g^k(X^{(i)})+ g(X^{(i)})]^2 &\leq 2\| g^k\|_n^2+2\| g\|_n^2
\leq 4 d(f)^2.
\end{eqnarray*}
Besides, the covering set is constructed such that 
$\| g^k-g\|_n \leq d(f)/8$. It follows that Step 1 is proved.

\medskip
\noindent
\underline{Proof of Step 2:}
We prove that for some constant $C$, 
$$\PX \left[ T_2 \geq \frac{d^2}{4}\right]\leq
\PX \left[\max_{1 \leq k \leq N_{\rm{pr}}} \{  d(f)^2-\|
  g^k\|_n^2  \}\geq \frac{d(f)^2}{4}\right]
 \leq
\exp\left(\log{N_{\rm{pr}}}- C  \frac{n d(f)^{2}}{1+d(f) +d(g)^{2}} \right).$$

As $g^{k} \in \cG'$, $d = \| g^k\|_2$. Then 
$$
 \max_{1\leq k\leq N_{\rm{pr}}} \left\{d(f)^2-\| g^k\|_n^2 \right\}=\max_{1\leq k\leq N_{\rm{pr}}}[\| g^k\|_2^2-\| g^k\|_n^2]
 .
 $$ 
Applying Theorem 3.5 in Chung and Lu \cite{ChungLu} we have
 that for all positive $\lambda$
 $$\PX \left[\sum_{i=1}^n   (g^k(\gX_{i}))^2\leq
    n\mathbb{E} (g^k(\gX_{i}))^2   -\lambda\right]\leq \exp\left(   -\frac{\lambda^2}{2n \mathbb{E}(g^k(\gX))^4}\right).$$
 Taking $\lambda=nd(f)^2/4$  and using that $\| g^{k}\|_2^2=d(f)^2$ we get
 $$\PX \left[\{  d^2-\| g^k\|_n^2  \}\geq \frac{d(f)^2}{4}\right]\leq \exp\Big(-\frac{nd(f)^4}{32\mathbb{E}(g^k(\gX))^4}\Big).$$
 It follows that 
 \begin{eqnarray}\label{step21}\PX \left[\max_{1\leq k\leq
       N_{\rm{pr}} } d(f)^2-\| g^k\|_n^2 
 \geq 
 \frac{d(f)^2}{4}    \right] \leq \sum_{k=1}^{N_{\rm{pr}}}
\exp\Big(-\frac{nd(f)^4}{32\mathbb{E}(g^k(\gX))^4}\Big) 
\leq \exp\Big(\log{N_{\rm{pr}}}-\frac{n d(f)^4}{32\max_{k}\mathbb{E}(g^k(\gX))^4}\Big).
 \end{eqnarray}

It remains to calculate $\EX g^{4}(\gX)$ for $g \in
\cG'$. Precisely we show the following result:
\begin{equation*}
 \mathbb{E}g^{4}(\gX)\leq 2 d(f)^{2}\left( 2 + 11 d(f)^{2} + 4
 d\right).
\end{equation*}
This result comes from the property of the RKHS $\cH$: indeed $g \in
\cH$ is written 
$g=\sum_{v  \in  \cP} g_{v}$ where the functions $g_{v}$ are centered
and orthogonal in $\bL^{2}(\PX)$. Therefore $\E g^{4}(\gX)$ is the
sum of the following terms:
\begin{eqnarray*}
A_1&=&\sum_{v\in \mathcal{P}} \EX g^4_v(\gX_v), \qquad
A_2 =\begin{pmatrix}4\\2\end{pmatrix} \sum_{v \not= v'}
\EX g^2_v(\gX_v)g^2_{v'}(\gX_{v'}),\\
A_3&=& 
\begin{pmatrix}4\\3\end{pmatrix} \sum_{v_1 \not= v_2\not= v_3}
\EX g^2_{v_1}(\gX_{v_1})
  g_{v_2}(\gX_{v_2}) g_{v_3}(\gX_{v_3}),
\quad A_4=\begin{pmatrix} 4\\3\end{pmatrix}
\sum_{v_1 \not= v_2}
\EX  g^3_{v_1}(\gX_{v_1})g_{v_2}(\gX_{v_2})\\
A_{5}&=&\begin{pmatrix}4\\1\end{pmatrix} \sum_{v_1 \not= v_2\not= v_3\not=v_4}
\EX g_{v_1}(\gX_{v_1})g_{v_2}(\gX_{v_2})g_{v_3}(\gX_{v_3})g_{v_4}(\gX_{v_4}).
\end{eqnarray*}
Using the Cauchy-Schwartz inequality and the fact that $\|
g_v\|_\infty\leq \| g_v\|_{\mathcal{H}_v}\leq 2$, and
$\|g\|_{2}=d(f)$ (because $g\in \cG'$), we get that $A_{1}$
is proportionnal to $d(f)^{2}$, $A_{2}, A_{3}, A_{5}$ to $d(f)^{4}$,
and $A_{4}$ to $d(f)^{3}$. For example 
\begin{equation*}
A_1=\sum_{v\in \mathcal{P}} \EX g^{4}_v(X_v)\leq \|
g\|_\infty^2 \sum_{v\in \mathcal{P}}
\|g_{v}\|_{2}^{2}= \|
g\|_\infty^2 
\|\sum_{v\in \mathcal{P}}g_{v}\|_{2}^{2} \leq 4 d(f)^2.
\end{equation*}
After calculation of the terms $A_i$, since $d(f)^{2}$ is assumed to be smaller than one,
we get that
\begin{eqnarray}
\label{step22}
\max_k \EX (g^k(\gX))^{4}\leq  4d(f)^2(1+d(f)
+(11/2) d(f)^{2}) \leq 34 d(f)^{2}.
\end{eqnarray}

Step 2 is proved  by combining \eref{step21} and \eref{step22}.

\medskip
\noindent
\underline{Proof of Step 3:} Let $N_{\rm{pr}}$ be defined at
Equation~\eref{CovN.eq}. We  prove that
$$\sqrt{\frac{\log N_{\rm{pr}}}{n}} \leq  \frac{64}{d(f)} \Ee \sup_{g\in \mathcal{G}'} \vert V_{n,\varepsilon}(g)\vert.$$
We start from \eref{encadreN} and write that
$$\log N_{\rm{pr}}(d(f)/8,\mathcal{G}',  
\|\cdot \|_{n} )\leq \log N(d(f)/16,\mathcal{G}',  
\|\cdot \|_{n} ).$$
Using the Sudakov minoration (see Pisier
\cite{Pisier1989}) we have
that for all positive $\omega$ 
$$\sqrt{\log N(\omega,\mathcal{G}', \|\cdot \|_{n} )} \leq
\frac{4 \sqrt{n}}{\omega}
\Ee \left[ \sup_{g\in \mathcal{G}'} \vert V_{n,\varepsilon}(g)\vert 
 \right].$$
Hence by taking $\omega=d(f)/16$, Step 3 is proved.\\

\noindent
\underline{Proof of Step 4:}\label{ProofS4}
The last step consists in bounding from above the Gaussian complexity  
$\mathbb{E}_\varepsilon  \sup_{g\in \mathcal{G}'} 
\vert V_{n,\varepsilon}(g)\vert$. 
This control is performed by using
Lemma  \ref{lemcomplex} (page~\pageref{lemcomplex}).
 According to Inequality \eref{but},
\begin{eqnarray*}
\vert V_{n,\varepsilon}(g)\vert \leq 
\kappa\left[ \sum_{v\in \mathcal{P}} \lambda_{n,v}^2 \| g_v\|_{\mathcal{H}_v}   + \sum_{v\in \mathcal{P}} \lambda_{n,v}  \| g_v\|_n  \right],
\end{eqnarray*}
with $\lambda_{n,v}$ defined by Equation~\eref{lambda.eq} satisfying 
$C_{1} \lambda_{n,v} \leq \gamma_{v}$ and $C_{1}\lambda_{n,v}^{2}\leq
\mu_{v}$ for all $v\in \cP$. 
It follows
\begin{eqnarray*}
\sup_{g\in \mathcal{G}'} \sum_{v \in \cP}\vert V_{n,\varepsilon}(g_v)\vert 
&\leq &
\kappa\sup_{g\in \mathcal{G}'} \left[ \sum_{v\in \mathcal{P}}\lambda_{n,v}^2 \| g_v\|_{\mathcal{H}_v}   +
\sum_{v\in \mathcal{P}} \lambda_{n,v}  \| g_v\|_n  \right]\\
&\leq & \frac{\kappa}{C_{1}}
\sup_{g\in \mathcal{G}'} \left[ \sum_{v\in \mathcal{P}}\mu_v\| g_v\|_{\mathcal{H}_v}   +
\sum_{v\in \mathcal{P}} \gamma_{v} \| g_v\|_n  \right]
\end{eqnarray*}

and  according to Condition \textbf{C1}, 
\begin{eqnarray*}
\sup_{g\in \mathcal{G}'} \sum_{v \in \cP}
 \left|V_{n,\varepsilon}(g_v) \right|
 &\leq&
\frac{4 \kappa}{C_{1}}
\left(\sup_{g\in \mathcal{G}'} \sum_{v \in{\mathcal S}}\mu_{v} \|g_{v}\|_{{\mathcal H}_{v}}
+  \sup_{g\in \mathcal{G}'} \sum_{v \in{\mathcal S}}\gamma_{v}\| g_v\|_n\right),
\\
 &\leq&\frac{4 \kappa}{C_{1}}
\left(2 \sum_{v \in{\mathcal S}}\mu_{v} 
+  \sup_{g\in \mathcal{G}'} \sum_{v \in{\mathcal S}}\gamma_{v}\| g_v\|_n\right),
\end{eqnarray*}
because $\| g_v\|_{\mathcal{H}_v} \leq 2$.
Now, according to Condition \textbf{C2}, and using that $2ab \leq
a^{2} + b^{2}$, we get
\begin{eqnarray*}
\sup_{g\in \mathcal{G}'} \sum_{v \in \cP}\vert V_{n,\varepsilon}(g_v)\vert 
&\leq &  \frac{4\kappa}{C_{1}}  \left[ 2\sum_{v\in S} \mu_v   +
 2 \sum_{v\in S} \gamma_v ^2+ \sup_{g\in \mathcal{G}'}\sum_{v\in S} \| g_v\|_2^2 \right]\\
 &\leq &\frac{4\kappa}{C_{1}}   \left[ 2\sum_{v\in S} (\mu_v   + \gamma_v^2) + d(f)^2\right],
\end{eqnarray*}
the last inequality coming from the fact for all $g \in \cG'$,
$\|g\|^{2}_{2} = d(f)^{2} \geq \sum_{v\in S} \| g_v\|_2^2 $.

Finally, thanks to~\eref{dn.eq}, we get
\begin{equation*}
 \sup_{g\in \mathcal{G}'} \sum_{v \in \cP}\vert
 V_{n,\varepsilon}(g_v)\vert 
\leq \frac{20 \kappa}{C_{1}} d(f)^{2}.
\end{equation*}

\hfill $\Box$

\subsection{\label{ProofsIntLemm.st}Proofs of intermediate Lemmas.
}

\subsubsection{\label{Prooflemcomplex}Proof of Lemma  \ref{lemcomplex} (page~\pageref{lemcomplex}): }
Let us write that the kernel $k_v$ is written as :
\begin{equation*}
 k_{v} (\gx_{v}, \gy_{v}) = \sum_{k} \omega_{v,k} 
\phi_{v,k}(\gx_{v})\phi_{v,k}(\gy_{v})
\end{equation*}
where $(\phi_{v,k})_{k=1}^{\infty}$
is an orthonormal basis of $\mathbb{L}^2(P_{v})$, where $P_v
= \prod_{a \in v} P_{a}$.

Let us consider the class of functions $\cK(t)$ defined as 
\begin{equation*}
  \cK(t) =\left\{ g_v \in \cH_{v}, \|g_v\|_{\cH_v}\leq 2, \|g_v\|_{2}\leq t \right\}.
\end{equation*}
It comes that
\begin{equation*}
g_{v} = \sum_{i} a_{i} \phi_{v,i}, 
\:\: \mbox{ with }
\|g_{v}\|^{2}_{\cH_{v}} = \sum_{i} \frac{a_{i}^{2}}{\omega_{v,i} }
\leq 4, \mbox{ and }
\|g_{v}\|^{2}_{2} = \sum_{i} a_{i}^{2} \leq t^{2}
\end{equation*}

In the following, we set  $\mu_{k,v} (t)= \min\left\{t^2, \omega_{k,v}\right\}$.  Hence
\begin{equation}
 \sum_{k} \frac{a_{k}^{2}}{\mu_{k,v}(t)} 
\leq 
  \frac{1}{t^{2}} \sum_{k}a_{k}^{2} + \sum_{k} \frac{a_{k}^{2}}{\omega_{k,v}}
= \frac{1}{t^{2}} \|g_{v}\|_2^{2} + 
\|g_{v}\|_{{\mathcal H}_{v}}^{2} \leq 5,
\label{ineq36.eq}
\end{equation}
as soon
  as  $g_{v}\in \cK(t)$.

Now, let us prove the lemma:
\begin{eqnarray*}
\EXe W_{n,2,v}(t)&=&
\EXe \sup_{g \in \cK(t) }
\left| \frac{1}{n}
\sum_{i=1}^n\varepsilon_i \sum_{\ell}a_{\ell}\phi_{v,\ell}(x_{v,i}) 
\right| \\
&=& \EXe  \sup_{g \in \cK(t)}
\left| 
\frac{1}{n}\sum_{\ell}\frac{a_{\ell}}{\sqrt{\mu_{v,\ell}(t)}} \sum_{i=1}^n\varepsilon_i
  \sqrt{\mu_{v,\ell} (t)}\phi_{v,\ell}(x_{v,i})\right| \\
& \leq & \sqrt{5}\sqrt{
\EXe  \sum_{\ell} \left( \frac{1}{n}\sum_{i=1}^n\varepsilon_i
  \sqrt{\mu_{v,\ell}(t) }\phi_{v,l}(x_{v,i})\right)^{2}} .
  \end{eqnarray*}
The last inequality follows from the Cauchy-Schwartz inequality and Inequality~\eref{ineq36.eq}.
  Now, simple calculation leads to 
\begin{equation*}
 \EXe  W_{n,2,v}(t) \leq \sqrt{5} \sqrt{\frac{1}{n}\sum_{\ell} \mu_{v,\ell} (t) },
\end{equation*}
\hfill $\Box$

\subsubsection{\label{Prooflemcompnormes1}Proof of Lemma
  \ref{lemcompnormes1} (page~\pageref{lemcompnormes1}): }
Using that $\left\vert \sqrt{a}-\sqrt{b}  \right\vert  \leq \sqrt{\vert a-b\vert},$
we get
$$\left\vert\| g_v\|_2- \| g_v\|_n\right\vert   \leq \sqrt{\left\vert \| g_v\|_2^2- \| g_v\|_n^2\right\vert} .$$
Hence
$$\left\lbrace  \| g_v\|_\infty\leq b , \;
\left\vert\| g_v\|_2- \| g_v\|_n\right\vert \geq \frac{b t}{2}\right\rbrace
\subset \left\lbrace \left\vert\| g_v\|_2^2- \| g_v\|_n^2\right\vert \geq \frac{b^2 t^2}{4}\right\rbrace.$$
The centered process
\begin{eqnarray*}
\left\vert\| g_v\|_2^2- \| g_v\|_n^2\right\vert&=&\left\vert\frac1n \sum_{i=1}^n g_v^2(\gX_{v,i})-\mathbb{E}(g_v^2(\gX_v))\right\vert,
\end{eqnarray*}
satisfies a concentration inequality given, for example,  by  Theorem
2.1 in Bartlett \textit{et al.} \cite{Bartlettetal2005} : if $\cC$ is a class of
functions $f$ such that $\|f\|_{\infty} \leq B$ and $E f(\gX)=0$, and
if there exists $\gamma>0$ 
such that for every $f \in \cC$, $\Var f(\gX) \leq \gamma^{2}$. Then for every
$x>0$, with probability at least $1-e^{-x}$, 
\begin{equation}
 \sup_{f\in \mathcal{C}}\frac{1}{n}\left\vert \sum_{j=1}^n
  f(\gX_j)\right\vert \leq \inf_{\alpha>0}\left\{
2(1+\alpha) E \left(\sup_{f\in \mathcal{C}}\frac{1}{n}\left| \sum_{j=1}^n
  f(\gX_j)\right| \right)
+ \sqrt{\frac{2  x}{n}} \gamma + B \left(\frac{1}{3} +
  \frac{1}{\alpha}\right)\frac{x}{n}
\right\}.
\label{Bartlet.eq}
\end{equation}

For any $t>0$, for $\cG(t)$ defined by~\eref{calG.eq},
let us consider  the class of functions $\cC(t)$
defined as follows
\begin{eqnarray*}
 \mathcal{C}(t) &=&\left\{f \mbox{ such that } f=g_v^2-\mathbb{E}(g_v^2), \mbox{
  with } g_v\in \cG(t) \right\}. \end{eqnarray*}
  Note that if $f \in \mathcal{C}(t)$, $\EX f(\gX_v)=0$ and $\|f\|_{\infty}
\leq  b^{2}$. We have to study 
\begin{eqnarray*}
 \gamma^{2}(t)  = \sup_{g_v \in \cG(t)} 
\EX \left(   g_v^2(\gX)-\|g_v\|^{2}_{2})\right)^{2} \mbox{ and }
\Gamma(t)  =  \EX \left( \sup_{g_v \in \cG(t)} \left|
    \|g_v\|^{2}_{n}-\|g_v\|^{2}_{2}\right| \right).
\end{eqnarray*}
It is easy to see that \begin{eqnarray*}
 \gamma^{2} (t)
& \leq &  b^{2} \sup_{g_v \in \cG(t)} \EX \left(
  g_v(\gX)+\|g_v\|_{2}\right)^{2} \leq 4 b^{2}t^{2}.
\end{eqnarray*}

Let  $\zeta_i$ be independent and identically  random variables
Rademacher distributed and let $E_{\gX,\zeta}$ denotes the
expectation with respect to the law 
of $(\gX,\zeta)$.
By  a symmetrization argument, 
$$\Gamma(t)  \leq 2 E_{\gX,\zeta}
\sup_{g_v\in \mathcal{G}(t)}
\left\vert\frac{1}{n} \sum_{i=1}^n \zeta_i g_v^2(\gX_i)\right\vert .$$

Since $\| g_v\|_{\infty}\leq b$, applying the  contraction principal (see
Ledoux-Talagrand \cite{ledoux1991}) we get that, for $Q_{n,v}(t)$ defined by
\eref{Qn},
\begin{eqnarray*}
E_{\gX,\zeta}
\sup_{g_v\in \mathcal{G}(t)}
\left\vert\frac{1}{n} \sum_{i=1}^n \zeta_i g_v^2(\gX_i)\right\vert
&\leq &4 b E_{\gX,\zeta}
\sup_{g_v\in \mathcal{G}(t)}
\left\vert\frac{1}{n} \sum_{i=1}^n \zeta_i g_v(\gX_i)\right\vert \\
&\leq & 4 b Q_{n,v}(t).
\end{eqnarray*}
The last inequality was 
proved by Mendelson~\cite{Mendelson2002}, Theorem 41 (see the proof of
Lemma~\ref{lemcomplex}).  Now, thanks
to~\eref{Bartlet.eq} we get that 
for all $x>0$, with probability greater than $1-e^{-x}$
\begin{equation*}
 \sup_{g_v \in \cG(t)}\left\vert \|g_{v}\|_{n}^{2} - \|g_{v}\|_{2}^{2}
  \right\vert \leq \inf_{\alpha>0}\left\{
16 (1+\alpha)  b Q_{n,v}(t)
+ \sqrt{\frac{2  x}{n}}2 b t + b^{2} \left(\frac{1}{3} +
  \frac{1}{\alpha}\right)\frac{x}{n}
\right\}.
\end{equation*}
Taking $x=c_{2} n t^{2}$, $t \geq  \nu_{n}$, we have that with
probability greater than $1-e^{-c_{2} n t^{2}}$
\begin{equation*}
 \sup_{g_v \in \cG(t)}\left\vert \|g_{v}\|_{n}^{2} - \|g_{v}\|_{2}^{2}
  \right\vert \leq \inf_{\alpha>0}t^{2} \left\{
16(1+\alpha) b \Delta 
+ \sqrt{2  c_{2}}4 b  + b^{2} \left(\frac{1}{3} +
  \frac{1}{\alpha}\right)c_{2}  
\right\}.
\end{equation*}

The infimum of the right hand side is reached in $\alpha=\sqrt{c_{2}
  b/ 16 \Delta}$,
and equals  
\begin{equation*}
 \frac{b^{2} c_{2}}{3} + 8 \sqrt{\Delta c_{2}} b^{3/2} + 4(4 \Delta + 
 \sqrt{2  c_{2}}) b.
\end{equation*}
The constants $\Delta$ and $c_{2}$ should satisfy that this infimum is
strictly smaller than $b^{2}/4$. 
For example, if $16 \Delta <b/8$, it remains to choose
$c_{2}$ small enough such that 
\begin{equation*}
b\left(\frac{c_{2}}{3} + \frac{\sqrt{2 c_{2}}}{2}\right) + 4\sqrt{2 c_{2}}
 < \frac{b}{8}.
\end{equation*}
\hfill \textbf{$\Box$}

\subsubsection{\label{Prooflemcompnormes2}Proof of Lemma \ref{lemcompnormes2} (page~\pageref{lemcompnormes2}): }
Let $t> \nu_{n,v}$ and $h$ be defined as $h = t g_v / \|g_v\|_{2}$. If
$g_v$ satisfies the assumptions of the lemma, then $h$  satisfies $\| h\|_2=t$, $\|h\|_{\mathcal{H}}\leq 2$ and $\|h\|_{\infty} \leq b$.
Applying Lemma \ref{lemcompnormes1}
(page~\pageref{lemcompnormes1}) to the function  $h$, we
obtain that for all $t \geq \nu_{n,v}$, with probability
greater than $1 - \exp(-c_2 n t^2)$, we have
$ |t - \|h\|_{n}| \leq bt/2$ for all $h \in \cG(t)$. This concludes the proof of the lemma.

\subsubsection{\label{Proofconcentration1}Proof of Lemma
  \ref{concentration1} (page~\pageref{concentration1}):}
The proof of  Lemma \ref{concentration1} is based on an isoperimetric
inequality for Gaussian processes (Borell \cite{Borell1975} or
Cirel'son \textit{et al.} \cite{cirelsonetal1976}) as it is  stated in
Theorem  (3.8), page 61 in Massart \cite{Massart}. Let us recall this inequality:

\begin{lem}
\label{lipschitz}
Let $P$ be the Gaussian probability measure on 
$\mathbb{R}^n$ and let $\phi$ be a function from $\mathbb{R}^n$ to
$\mathbb{R}$, and $\|\phi\|_{L}$ its Lipschitz semi-norm:
$$\|\phi\|_{L}=\sup_{x\not=y}\frac{\vert \phi(x)-\phi(y)\vert}{\sqrt{n}\|
  x-y\|_{n}}.$$  Let $\overline{\Phi}$ be the cumulative distribution of the standard Gaussian distribution.
Then for any $u$,
\begin{eqnarray}
P(\vert f-E_P f \vert \geq u\ )\leq 4 \overline{\Phi}\left( \frac{u}{\|\phi\|_{L}}   \right).
\end{eqnarray}
\end{lem}
We apply Lemma~\ref{lipschitz} to $\phi(\varepsilon_{1}, \ldots,
\varepsilon_{n}) = W_{n,n,v}(t)$. By Cauchy-Schwarz Inequality, 
$\|\phi\|_{L}=t/\sqrt{n}$.
It follows that  Lemma \ref{concentration1} is proved since
\begin{eqnarray*}
\PXe \left(\vert W_{n,n,v}( t)-\mathbb{E}_{\varepsilon} W_{n,n,v}( t)\vert \geq \delta t\right) \leq 4\exp\left\lbrace
 -\frac{(\delta t)^2 }{2\left(\frac{t}{\sqrt{n}} \right)^2}\right\rbrace \leq 4\exp\left( -\frac{n\delta^2}{2}  \right).
 \end{eqnarray*}
\hfill $\Box$

\subsubsection{\label{Proofconcentration2}Proof of Lemma \ref{concentration2} (page~\pageref{concentration2}):}
We start with the proof of \eref{concentn2} in Lemma
\ref{concentration2} by applying once again Lemma~\ref{lipschitz}
given above, to
the function $\phi(\geps) = \phi(\varepsilon_{1}, \ldots,
\varepsilon_{n}) = W_{n,2,v}(t)$.
On the event  $\Omega_{v, t}$ defined by~\eref{Omega},
we have $\|g_{v}\|_{n} \leq bt/2 + \|g_{v}\|_{2}$. Besides if
$\|g_{v}\|_{\cH_{v}} \leq 2$, then $\|g_{v}\|_{\infty} \leq
2$. Therefore applying Lemma~\ref{lemcompnormes1} with $b=2$, we get
that if $\|g_{v}\|_{2} \leq t$,
\begin{eqnarray*}
\left| \phi(\geps) - \phi(\geps')\right| \leq \sup_{\|g_v\|_{n} \leq 2t}\|g_v\|_{n} \|\geps-\geps' \|_{n}
\leq  2 t \|\geps-\geps' \|_{n},
\end{eqnarray*} 
leading to $\|\phi\|_{L}=2 t/ \sqrt{n}$. 
It follows that  \eref{concentn2} in Lemma  \ref{concentration2} is proved since
$$\PXe \left[\left\lbrace\vert W_{n,2,v}( t)-
\mathbb{E}_{\varepsilon}W_{n,2,v}( t)
\vert    \geq \delta t\right\rbrace \cap \Omega_{v,t}^c\right] \leq 4\exp\left\lbrace
 -\frac{(\delta t)^2 }{2\left(\frac{2t}{\sqrt{n}} \right)^2}\right\rbrace \leq 4\exp\left( -\frac{n\delta^2}{8}  \right).$$

\medskip
\noindent
We now come to the proof of \eref{concentration} in Lemma
\ref{concentration2} 
using a Poissonian inequality for self-bounded processes (see
Boucheron \textit{et al.} \cite{Boucheronetal2000}) and  Theorem 5.6, p 158 in
Massart \cite{Massart}). Let us recall it in the particular
case we are interested in:

\begin{theo}
\label{autoborne}
Let  $X_1,\cdots, X_n$  be $n$ iid variables. 
For
$i\in \{1,\cdots,n\}$ let  $X_{(-i)}=(X_1,\ldots, X_{i-1},X_{i+1},\ldots,X_n)$.
Let $Z$ be a nonnegative and bounded measurable function of $X=(X_1,\cdots,X_n)$. Assume that for all $i\in\{1,\cdots,n\}$, there exists a measurable function
$Z_i$ of  $X_{ (-i)}$ such that
$0<Z-Z_i\leq 1,$ and 
 $\sum_{i=1}^n (Z-Z_i)\leq Z$.
 Then, for 
all  $x>0$, we have $P \{Z\geq E(Z) +x \}  \leq \exp\left(
  -x^2/2 E(Z)   \right)$.
\end{theo}

We apply this result to $Z$ defined as 
\begin{equation*}
Z=Z(\gX_1,\cdots,\gX_n)= 
n \Ee W_{n,2,v}(t) = 
n \Ee
\sup \left\{\vert V_{n,\varepsilon}(g_v)\vert, \| g_v\|_2\leq t, 
\| g_v\|_{\cH_v}\leq 2\right\}   .
\end{equation*}
The variable $Z$ is positive, and because the distribution of $(\varepsilon_{1}, \ldots,
\varepsilon_{n})$ is symmetric, we have that 
\begin{equation*}
Z=\Ee 
\sup \left\{ n V_{n,\varepsilon}(g_v), \| g_v\|_2\leq t, 
\| g_v\|_{\cH_v}\leq 2\right\}.
\end{equation*}

Let $\tau$ be the function in $\mathcal{H}_v$ such that 
$Z= \Ee n V_{n,\varepsilon}(\tau)$ (note that $\tau$
depends on $(\gX_1, \ldots, \gX_n)$ and on $(\varepsilon_{1}, \ldots,
\varepsilon_{n})$), and let 
\begin{equation*}
Z_{i}=  \Ee  \sup_{g_v} 
\sum_{j\not =i} \varepsilon_j g_{v}(\gX_j).
\end{equation*}
 We  show that $Z$ and $Z_{i}$ satisfy the assumptions of Theorem~\ref{autoborne}: 
\begin{eqnarray*}
 Z-Z_{i} &=  &
\Ee 
\left(
\varepsilon_i \tau( \gX_i) +
\sum_{j \neq i}\varepsilon_j \tau( \gX_j)
- \sup_{g_{v}} \sum_{j \neq i} \varepsilon_jg_{v}( \gX_j) \right) \\
&\leq  & \Ee  \varepsilon_i \tau( \gX_i) \leq
\frac{1}{\sqrt{2\pi}} \Ee \sup_{x \in \cX}|\tau(\gX)|
\leq \sqrt{\frac{2}{\pi}},
\end{eqnarray*}
where the last inequality comes from the fact that $\sup_{x \in
  \cX}|\tau(\gX)| \leq \|\tau\|_{\cH_v} \leq 2$.
Moreover  $Z-Z_i \geq 0$ since 
\begin{eqnarray*}
Z & = & 
 \Ee  
\sup_{g_{v}} \sum_{j=1}^{n} \varepsilon_j g_v( \gX_j) = 
 \Ee  \left(E_{\varepsilon_{i}} \sup_{g_{v}}
\sum_{j=1}^{n} \varepsilon_j g_v( \gX_j) \right)\\
&\geq&   \Ee  \left(\sup_{g_{v}}
E_{\varepsilon_{i}}\sum_{j=1}^{n} \varepsilon_j g_v( \gX_j) \right)= Z_{i}.
\end{eqnarray*}

Finally we have:
$$\sum_i (Z-Z_{i})= \sum_{i=1}^n \Ee   
\left( \varepsilon_i\tau(X_i)+ \sum_{j\not=i}^n \varepsilon_j
  \tau(X_j) 
-\sup_{g-v}\sum_{j\not =i}^n \varepsilon_j g_v(X_j) \right)
\leq \sum_{i=1}^n \Ee  \varepsilon_i\tau(X_i) = Z.$$

Therefore, following Theorem~\ref{autoborne}, we get that for all
postive $u$
\begin{eqnarray*}
\PXe \left[ \Ee  W_{n,2,v}(t)
    -\EXe W_{n,2,v}(t) \leq  \frac{u}{n} \right]\leq \exp\left[  -\frac{u^2}{
\EXe W_{n,2,v}(t)}   \right].
\end{eqnarray*}
As $\EXe  W_{n,2,v}(t)\leq Q_{n,v}(t)$, see Lemma~\ref{lemcomplex} page~\pageref{lemcomplex}, we get the expected result since for all
positive $x$
\begin{eqnarray*}
\PX \left[\mathbb{E}_\varepsilon W_{n,2,v}(t)\geq \EXe W_{n,2,v}(t)+x\right]\leq \exp\left[  -\frac{n x^2}{
Q_{n,v}(t)}   \right].
\end{eqnarray*}

\hfill $\Box$

\subsubsection{\label{ProofCase1}Proof of Lemma \ref{Case1} (page~\pageref{Case1}):}   
From Lemma \ref{concentration1}, page~\pageref{concentration1} with $t=\lambda_{n,v}=\delta$, we get that with  probability
 greater than $1-4\exp(-n\lambda_{n,v}^2/2)$,
$$ \Ee W_{n,n,v}( \lambda_{n,v})\leq \lambda_{n,v}^2+W_{n,n,v}( \lambda_{n,v} ).$$
 
The next step consists in comparing $W_{n,n,v}( \lambda_{n,v})$ and $
W_{n,2,v}( 2\lambda_{n,v})$. Recall that $\lambda_{n,v} \geq
\nu_{n,v}$, see~\eref{lambda.eq}. Let $g_{v}$ such that $\|g_v\|_{n}
\leq \lambda_{n,v}$.
\begin{itemize}
\item 
When $\| g_v\|_2 \leq \lambda_{n,v}$,  according to Lemma \ref{lemcompnormes1}
(page~\pageref{lemcompnormes1}), taking $b=2$ , since since $\| g_v\|_n\leq \lambda_{n,v}$, we get that 
with probability greater than $1-\exp(-c_{2} n \lambda_{n,v}^{2})$, 
$$\| g_v\|_n -\lambda_{n,v}\leq \| g_v\|_2 \leq
\|g_v\|_n+\lambda_{n,v} \leq 2 \lambda_{n,v}.$$

\item When $\| g_v\|_2 \geq t$, we apply Lemma
\ref{lemcompnormes2} (page~\pageref{lemcompnormes2}) with $b=2$.
For any function $g_{v}$ such that
$\| g_{v}\|_\infty \leq 2$, and $\| g_{v}\|_2\geq \lambda_{n,v}$, we
have $\| g_{v}\|_2\leq 2 \| g_{v}\|_n \leq 2 \lambda_{n,v}$.
\end{itemize}
This implies that,  with probability greater than $1-\exp(-c_2n \lambda_{n,v}^2)$
 we have
$$W_{n,n,v}( \lambda_{n,v}) \leq W_{n,2,v}( 2\lambda_{n,v}).$$

We now study the process $W_{n,2,v}( \lambda_{n,v})$. By applying
\eref{concentn2} in Lemma  \ref{concentration2},
page~\pageref{concentration2}, with $\delta=t=\lambda_{n,v}$ we get
that with  probability greater than $1-4\exp(- n\lambda_{n,v}^2/8)$
$$ W_{n,2,v}( \lambda_{n,v})\leq \lambda_{n,v}^2+ \Ee (W_{n,2,v}( \lambda_{n,v})).$$
It follows that
\begin{eqnarray*}
\Ee W_{n,n,v}( \lambda_{n,v})&\leq& \lambda_{n,v}^2+W_{n,n,v}( \lambda_{n,v} ))\\
&\leq&\lambda_{n,v}^2+W_{n,2,v}( 2\lambda_{n,v}))\\
&\leq& 5\lambda_{n,v}^2+\Ee  (W_{n,2,v}( 2\lambda_{n,v})).
\end{eqnarray*}
Next, we apply \eref{concentration} in Lemma \ref{concentration2}, with
$t=2\lambda_{n,v}$ and $x= 4 \lambda^{2}_{n,v}$. We get that $$
\Ee W_{n,2,v}( 2\lambda_{n,v})\leq
4\lambda_{n,v}^2+ \EXe  (W_{n,2,v}(
2\lambda_{n,v})),$$
with probability greater than 
\begin{equation*}
 1 -2 \exp\left(-16 \frac{n \lambda^{4}_{n,v}}{Q_{n,v}(2
     \lambda_{n,v})}\right)
\geq 1 -2 \exp\left(- \frac{4 n \lambda^{2}_{n,v}}{\Delta }\right).
\end{equation*}
The last inequality comes from the definition of $\nu_{n,v}$,
see~\eref{nu}, and from the fact that $\lambda_{n,v} \geq \nu_{n,v}$,
see~\eref{lambda.eq}.

Putting everything together, we get that with probability greater than
$1 - c_{1}\exp(-c_2 n\lambda_{n,v}^2)$ for some positive constants
$c_1, c_2$, 
\begin{eqnarray*}
\Ee W_{n,n,v}( \lambda_{n,v})
&\leq& 9\lambda_{n,v}^2+ \EXe  (W_{n,2,v}( 2\lambda_{n,v}))\\
&\leq& 9\lambda_{n,v}^2+ Q_{n,v}(2\lambda_{n,v}),
\mbox{ thanks to Lemma~\ref{lemcomplex}, page~\pageref{lemcomplex}},
\\
&\leq& 9\lambda_{n,v}^2+ 4 \Delta \lambda_{n,v}^2.
\end{eqnarray*}
Applying once again  Lemma \ref{concentration1},
page~\pageref{concentration1}, we get that
\begin{equation*}
W_{n,n,v}(\lambda_{n,v}) \leq   \Ee  W_{n,n,v}(
\lambda_{n,v})+ \lambda_{n,v}^2 \leq   \Big(10+4 \Delta\Big)\lambda_{n,v}^2. 
\end{equation*}
This ends the proof of the lemma by taking $\kappa=10+4 \Delta$.

\hfill \textbf{$\Box$}

\subsection{\label{algoProofs.st}Algorithm:
  Propositions~\ref{Casemu0.pr}-\ref{Solvteta.pr}}
We consider the minimization of $C'(f_{0}, \gtheta)$ given at
Equation~\eref{CritAlgo.eq}. Because $C'(f_{0}, \gtheta)$ is convex and separable, we use a block
coordinate descent algorithm, group $v$ by group $v$. We refer to
Bubeck~\cite{Bubeck2014} for a review on convex
optimization.

In what follows, the group  $v$ is
fixed, and for given values of
$f_{0}$ and $\gtheta_w, w \neq v$, we look for the minimizer of $C'$
with respect to $\gtheta_v$: Setting
\begin{equation*}
 \gR_{v}=\gY - f_{0}-\sum_{w\neq v}K_{w}\gtheta_{w},
\end{equation*}
we aim at minimizing with respect to $\gtheta_{v}$
\begin{equation*}
Q( \gtheta_v) = \|\gR_{v} - K_{v}\gtheta_v\|^{2} + \gamma^{'}_{v}
\|K_{v} \gtheta_v\| + \mu^{'}_{v} \|K^{1/2}_{v} \gtheta_v\|.
\end{equation*}

If $\gamma'_{v}=\mu'_{v}=0$, then $\gtheta_v=K_{v}^{-1}\gR_{v}$ is
the solution. In what follows we consider the case where at least one
of both is non zero.

If $\partial Q_{v}$ denotes the subdifferential of
$Q(\gtheta_v )$ with respect to  $\gtheta_{v}$, we need to solve
$0 \in \partial Q_{v}$. Let us recall that for all $v \in \cP$ the
matrices $K_{v}$ are symmetric and  strictly definite positive.

Let us begin with the calculation of the subdifferential of 
$\|K_{v}\theta_{v}\|$ with respect to $\theta_{v}$: if $\theta_{v} \neq 0$,
we have
\begin{equation*}
 \frac{\partial \|K_{v}\theta_{v}\|}{\partial \theta_{v}} =
 \frac{K_{v}^{2}\theta_{v} }{\|K_{v}\theta_{v}\|},
\end{equation*}
and if $\theta_{v}=0$ the subdifferential is the set of $x \in
\R^{n}$ such that $\|K_{v}^{-1}x\| \leq 1$. 

Therefore if $\theta_{v}
\neq 0$, 
\begin{eqnarray*}
 \partial Q_{v} 
&= &-2 K_{v} \gR_{v} + 2 K_{v}^{2}\theta_{v}  + \gamma'_{v} 
\frac{K_{v}^{2}\theta_{v} }{\|K_{v}\theta_{v}\|}+ \mu'_{v} \frac{K_{v}\theta_{v} }{\|K_{v}^{1/2}\theta_{v}\|},
\end{eqnarray*} 
while if $\theta_{v} = 0$, 
\begin{equation*}
 \partial Q_{v} = \left\{   -2 K_{v}\gR_{v} + \gamma'_{v} t +
\mu'_{v} s \mbox{ where }  t, s \in \R^{n} \mbox{ such that }
\|K_{v}^{-1} t\| \leq 1, \|K_{v}^{-1/2} s\|\leq 1
\right\}.
\end{equation*}

Let us begin with the case $\mu'_{v}=0$. 
\begin{prop}
Let $( \rho)_{+}$ denotes the positive part of $\rho\in \R$. If $\mu'_v=0$, 
\begin{equation*}
\theta_{v} = \left(1 - \frac{\gamma'_v}{\|2 R_v\|}\right)_{+} 
K_{v}^{-1} \gR_{v}.
\end{equation*}
\label{Casemu0.pr}
\end{prop}

\paragraph{Proof of Proposition~\ref{Casemu0.pr}}
The problem comes to
  minimize 
\begin{equation*}
U( \gbeta_v) = \|\gR_{v} - \gbeta_v\|^{2} + \gamma^{'}_{v}
\|\gbeta_v\| 
\end{equation*}
and to take $\gtheta_{v} = K_{v}^{-1} \gbeta_v$. The subdifferential
of $U$ is given by 
\begin{eqnarray*}
 \partial U_{v}(\gbeta_{v} )
&= &-2  \gR_{v} + 2\gbeta_{v}  + \gamma'_{v} 
\frac{\gbeta_{v} }{\|\gbeta_{v}\|} \mbox{ if } \gbeta_{v}\neq 0 \\
\partial U_{v}(0) & = &  \left\{   -2 \gR_{v} + \gamma'_{v} t \mbox{ where }  
\| t\| \leq 1
\right\}.
\end{eqnarray*}
We get 
\begin{equation*}
 \gbeta_{v} = \left(1 - \frac{\gamma'_v}{\|2 R_v\|}\right)_{+} \gR_{v}.
\end{equation*}
\hfill \textbf{$\Box$}

Let $\mu'_{v} >0$, the following proposition gives a
necessary and sufficient condition for $\gtheta_{v} =0$ to be the
minimizer of $Q$. For the sake of clarity let us recall the
definitions of $J$ and $J^{*}$ given in Section~\ref{algo.st}: For all $x \in \R^{n}$
\begin{equation*}
J(x) = \| 2 \gR_{v}  - \mu'_{v} K_{v}^{-1} x\|^{2}, \mbox{ and }
J^{*} = \min \left\{ J(x), \mbox{ for } x \mbox{ such that } 
\| K_{v}^{-1/2} x\| \leq 1\right\}.
\end{equation*}

\begin{prop}
The minimizer of $Q(\gtheta_v)$ is  $\gtheta_{v}=0$ 
 if and only if $J^{*} \leq \gamma'^{2}_{v}$.
\label{CNS.pr}
\end{prop}

\paragraph{Proof of Proposition~\ref{CNS.pr}}

\begin{enumerate}
\item Let us assume that $\gtheta_{v}=0$ is 
the minimizer of $Q(\gtheta_v)$. Then $0 \in \partial Q$ implies that  there exists $t^{*}$ and $s^{*}$ such that $\|K_{v}^{-1} t^{*}\| \leq
1$ and $\|K_{v}^{-1/2} s^{*}\| \leq 1$ and such that $\gamma'_{v} t^{*}+
\mu'_{v} s^{*}= 2 K_{v}\gR_{v}$.

If $\gamma'_{v} > 0$, then 
\begin{equation*}
 K_{v}^{-1} t^{*}= \frac{1}{\gamma'_{v}} \left(2 \gR_{v} - \mu'_{v}
   K_{v}^{-1}s^{*}\right), \mbox{ and } \|K_{v}^{-1} t^{*}\| = \frac{1}{\gamma'_{v}} \sqrt{J(s^{*})}.
\end{equation*}
Because $J^{*} \leq J(s^{*})$ and $\|K_{v}^{-1} t^{*}\| \leq 1$, we get that $J^{*} \leq \gamma'^{2}_{v}$.

If $\gamma'_{v}=0$, then $\mu'_{v} s^{*}= 2 K_{v}\gR_{v}$ and
$J(s^{*})=J^{*}=0$.
\item Let us now assume that $J^{*} \leq \gamma'^{2}_{v}$. Note that
minimizing the convex function $J(s)$ over the convex set
$\|K_{v}^{-1/2}s\| \leq 1$ has always a solution. Let us denote this
solution by $s^{*}$. 

If $\gamma^{'}_{v} >0$, let $t^{*} = (2 K_v \gR_{v} - \mu'_{v}
s^{*})/\gamma^{'}_{v}$. Then $-2 K_{v}\gR_{v}+\gamma'_{v} t^{*}+
\mu'_{v} s^{*}= 0$, and $-2 K_{v}\gR_{v}+\gamma'_{v} t^{*}+
\mu'_{v} s^{*} \in \partial Q(0)$ since $\|K_{v}^{-1/2} s^{*}\| \leq 1$, and 
$\|K_{v}^{-1} t^{*}\| = J(s^{*})/\gamma^{'}_{v} \leq 1$. Therefore $\gtheta_v=0$
is the minimizer of $Q$. 

If $\gamma'_v=0$, and $J^{*}=0$, then  $-2 K_{v}\gR_{v}+
\mu'_{v} s^{*}= 0$, and  $\|K_{v}^{-1/2} s^{*}\| \leq 1$. Therefore $\gtheta_v=0$
is the minimizer of $Q$. 

\end{enumerate}
\hfill \textbf{$\Box$}

\begin{prop}
\label{RdgPb.pr}
Let $\mu'_v > 0$ and $\gtheta_{v}$ be the minimizer of $Q$. 
\begin{enumerate}
\item If $\|2 K_{v}^{1/2}\gR_{v}\| \leq \mu'_v$, then $\gtheta_{v}=0$
\item If not, for $\rho>0$, let 
\begin{equation*}
\gtheta(\rho)  = 2 \mu'_{v} \left(\mu'^{2}_{v} K_{v}^{-1} + \rho I_{n}\right)^{-1} K_v^{-1/2}R_v,
\end{equation*}
and let $\rho^{*}$ defined as $\|\gtheta(\rho^{*})\|=1$. 
Then $J(\theta(\rho^{*})) \leq \gamma'^{2}_v$ if and only if  $\gtheta_{v}=0$, 
\end{enumerate}
\end{prop}

\paragraph{Proof of Proposition~\ref{RdgPb.pr}}
Minimizing $J(x)$ under the constraint $\|K_{v}^{-1/2}x\| \leq 1$ is
equivalent to minimizing $K(\gbeta)=\|2 \gR_{v} - \mu'_v K_{v}^{-1/2}\gbeta\|^{2}$
under the constraint $\|\gbeta\|^{2}\leq 1$. Let $\gbeta^{*}=2
K_{v}^{1/2}\gR_{v}/\mu'_v$. Then $K(\gbeta^{*}) =0$, which is smaller
than 
$\gamma'^{2}_{v}$,  and if $\|\gbeta^{*}\| \leq 1$, following
Proposition~\ref{CNS.pr}, we get $\gtheta_{v}=0$.

If $\|\gbeta^{*}\| > 1$, we have to solve a ridge regression problem
by minimizing $K(\gbeta) + \rho \|\gbeta\|^{2}$ for some positive
$\rho$. The solution is given by $\gtheta(\rho)$. Let us note that
$\|\gtheta(\rho)\|$ decreases to 0 when $\rho$ tends to infinity and
that its maximum is $\|\gtheta(0)\| = \|\gbeta^{*}\|$. Therefore if
$\|\gbeta^{*}\|>1$, there exists $\rho^{*}$ such that
$\|\gtheta(\rho^{*})\|=1$. Following Proposition~\ref{CNS.pr},
Proposition~\ref{RdgPb.pr}  is proved. A numerical procedure can be
used for calculating this $\rho^{*}$. \hfill \textbf{$\Box$}

\medskip

Let us now consider the case where $\gtheta_{v}$ is non zero. It should satisfy the subgradient
condition $\partial Q_{v} = 0$ which leads to
\begin{equation*}
\theta_{v} = \left(2 K_{v}^{2}  + \gamma'_{v} 
\frac{K_{v}^{2} }{\|K_{v}\theta_{v}\|}+ \mu'_{v}
\frac{K_{v} }{\|K_{v}^{1/2}\theta_{v}\|}\right)^{-1}
2 K_{v} \gR_{v}
\end{equation*}
which is equivalent to Equation~\eref{subdiffNzero.eq}.

\begin{prop}
\label{Solvteta.pr}
For all $\rho_{1}, \rho_{2}>0$ let 
\begin{equation*}
 \gtheta(\rho_{1}, \rho_{2}) = \left(K_{v} + \rho_{1} K_{v} + \rho_{2}
 I_{n}\right)^{-1} \gR_{v}.
\end{equation*}
If $\mu'_{v} >0$,  there exists a non zero solution to
Equation~\eref{subdiffNzero.eq} if and only if there exists $\rho_{1},
\rho_{2}>0$ such that 
\begin{equation}
\left.  
\begin{array}{rcl}
\gamma'_{v} & = & 2 \rho_{1} \|K_{v} \gtheta(\rho_{1}, \rho_{2})\| \\
\mu'_{v} & = & 2 \rho_{2} \|K^{1/2}_{v} \gtheta(\rho_{1}, \rho_{2})\| 
\end{array}
\right\} \label{Solvalpha.eq}
\end{equation}
Then $\gtheta_{v} = \gtheta(\rho_{1}, \rho_{2})$.
\end{prop}

\paragraph{Proof of Proposition~\ref{Solvteta.pr}}
If there exists a non zero solution to
Equation~\eref{subdiffNzero.eq} then $\|K_{v} \gtheta_{v}\|$ and
$\|K_{v}^{1/2} \gtheta_{v}\|$ are non zero because $K_{v}$ is definite
positive. 
Let $\rho_{1} = \gamma'_{v}/2\|K_{v} \gtheta_{v}\|$ and $\rho_{2} =
\mu'_{v}/2\|K^{1/2}_{v} \gtheta_{v}\|$, then  
\begin{equation*}
 \gtheta(\rho_{1}, \rho_{2}) = \left(K_{v} + 
\frac{\gamma'_{v}}{2\|K_{v} \gtheta_{v}\|} K_{v} + 
\frac{\mu'_{v}}{2\|K^{1/2}_{v} \gtheta_{v}\|}
 I_{n}\right)^{-1} \gR_{v} = \gtheta_{v},
\end{equation*}
and, for such $\rho_{1}, \rho_{2}$, Equation~\eref{Solvalpha.eq}
is satisfied.

Conversely, if there exist $\rho_{1}, \rho_{2}$ such that
Equation~\eref{Solvalpha.eq} is satisfied, then necessarily $\|K_{v}
\gtheta(\rho_{1}, \rho_{2)})\|$ and
$\|K_{v}^{1/2} \gtheta(\rho_{1}, \rho_{2)})\|$ are non zero and 
$\rho_{1}=\gamma'_{v}/2 \|K_{v}\gtheta(\rho_{1}, \rho_{2})\|$ and 
$\rho_{2}=\mu'_{v}/2 \|K_{v}^{1/2}\gtheta(\rho_{1}, \rho_{2})\|$. Then
\begin{equation*}
 \gtheta(\rho_{1}, \rho_{2}) = 
\left(K_{v} 
+ \frac{\gamma'_{v}}{2\|K_{v}\gtheta(\rho_{1}, \rho_{2})\|}K_{v} 
+ \frac{\mu'_{v}}{2\|K_{v}^{1/2}\gtheta(\rho_{1}, \rho_{2})\|}
 I_{n}\right)^{-1} \gR_{v},
\end{equation*}
which is exactly Equation~\eref{subdiffNzero.eq} calculated in
$\gtheta_{v} = \gtheta(\rho_{1}, \rho_{2})$. \hfill \textbf{$\Box$}
 
\medskip

\paragraph{{\bf Taking into account  that $\|K_{v^{1/2}} \gtheta_{v}\|
   \leq r_{v}$}}

As already mentionned in Section~\ref{CalcEstim.st}, one may want to
minimize $C(f_{0}, \gtheta)$ under the additional constraint that
$\|K_{v^{1/2}} \gtheta_{v}\| \leq r_{v}$ for some positive constant
$r_{v}, v \in \cP$. 

For each group $v$, we have thus to minimize $Q(\gtheta_{v})$ under the
constraint $\|K_{v^{1/2}} \gtheta_{v}\| \leq r_{v}$. We know that this
problem is equivalent to minimize 
\begin{equation*}
 Q(\gtheta_{v}) +
\lambda\|K_{v}^{1/2} \gtheta_{v}\| =
\|\gR_{v} - K_{v}\gtheta_v\|^{2} + \gamma^{'}_{v}
\|K_{v} \gtheta_v\| + (\mu^{'}_{v} + \lambda) \|K^{1/2}_{v} \gtheta_v\|.
\end{equation*}
for some $\lambda$ that depends on $r_v$. 

Let us first remark that, for a fixed $\gamma^{'}_{v}$, and $\lambda
\geq 0$, if 
$\widehat{\gtheta}_{v}(\mu^{'}_{v} + \lambda)$ minimizes $Q(\gtheta_{v}) +
\lambda\|K_{v}^{1/2} \gtheta_{v}\|$ with respect to $\gtheta_{v}$,
then $\|K^{1/2}_{v} \widehat{\gtheta}_{v}(\mu^{'}_{v} + \lambda) \|
\leq \|K^{1/2}_{v} \widehat{\gtheta}_{v}(\mu^{'}_{v})\|$. It can be
easily proved by writing
\begin{eqnarray*}
 Q\left(\widehat{\gtheta}_{v}(\mu^{'}_{v} + \lambda)\right) +
  \lambda\|K^{1/2}_{v} \widehat{\gtheta}_{v}(\mu^{'}_{v} + \lambda) \|
  & \leq &
Q\left(\widehat{\gtheta}_{v}(\mu^{'}_{v})\right) +
  \lambda\|K^{1/2}_{v} \widehat{\gtheta}_{v}(\mu^{'}_{v}) \| \\
 & \leq &Q\left(\widehat{\gtheta}_{v}(\mu^{'}_{v} + \lambda)\right) +
\lambda\|K^{1/2}_{v} \widehat{\gtheta}_{v}(\mu^{'}_{v}) \|.
\end{eqnarray*}

Therefore, one can proceed as follows: calculate $\widehat{\gtheta}_{v}$ for
$\lambda=0$. If $\|K^{1/2}_{v} \widehat{\gtheta}_v\| \leq r_{v}$, then
one go to the next step of the algorithm. If $\|K^{1/2}_{v}
\widehat{\gtheta}_v\| >  r_{v}$, one increases $\lambda$ untill $\|K^{1/2}_{v}
\widehat{\gtheta}_v\| =  r_{v}$. 

\bibliography{biblioArticle}

\end{document}